\documentclass[12pt,a4paper]{article}

\usepackage[utf8]{inputenc}
\usepackage[german,french]{babel}
\usepackage[utf8]{inputenc}
\usepackage{pgf,tikz}
\usetikzlibrary{arrows}
\usepackage[cyr]{aeguill}

\usepackage{graphicx}

\newcommand{\la}{\langle} 
\newcommand{\ra}{\rangle}

\newcommand{\Bilin}{\mathop{\rm Bilin}\nolimits}
\newcommand{\CI}{\mathop{\rm CI}\nolimits}
\newcommand{\pol}{\mathop{\rm pol}\nolimits}
\newcommand{\vect}{\wedge}
\newtheorem{definition}{\noindent \textbf{D\'efinition}}[section]

\newcommand{\PP}{\mathop{\mathbf{P}}\nolimits}
\newcommand{\Hom}{\mathop{\rm Hom}\nolimits}
\newcommand{\cQ}{\mathop{\mathcal{Q}}\nolimits}
\newcommand{\cB}{\mathop{\mathcal{B}}\nolimits}

\renewcommand{\bar}{\overline}

\newcommand{\cC}{\mathcal{C}}
\newcommand{\cD}{\mathcal{D}}

\newcommand{\vs}{\vspace*{3mm}}

\DeclareMathSymbol{A}{\mathalpha}{operators}{`A}
\DeclareMathSymbol{B}{\mathalpha}{operators}{`B}
\DeclareMathSymbol{C}{\mathalpha}{operators}{`C}
\DeclareMathSymbol{D}{\mathalpha}{operators}{`D}
\DeclareMathSymbol{E}{\mathalpha}{operators}{`E}
\DeclareMathSymbol{F}{\mathalpha}{operators}{`F}
\DeclareMathSymbol{G}{\mathalpha}{operators}{`G}
\DeclareMathSymbol{H}{\mathalpha}{operators}{`H}
\DeclareMathSymbol{I}{\mathalpha}{operators}{`I}
\DeclareMathSymbol{J}{\mathalpha}{operators}{`J}
\DeclareMathSymbol{K}{\mathalpha}{operators}{`K}
\DeclareMathSymbol{L}{\mathalpha}{operators}{`L}
\DeclareMathSymbol{M}{\mathalpha}{operators}{`M}
\DeclareMathSymbol{N}{\mathalpha}{operators}{`N}
\DeclareMathSymbol{O}{\mathalpha}{operators}{`O}
\DeclareMathSymbol{P}{\mathalpha}{operators}{`P}
\DeclareMathSymbol{Q}{\mathalpha}{operators}{`Q}
\DeclareMathSymbol{R}{\mathalpha}{operators}{`R}
\DeclareMathSymbol{S}{\mathalpha}{operators}{`S}
\DeclareMathSymbol{T}{\mathalpha}{operators}{`T}
\DeclareMathSymbol{U}{\mathalpha}{operators}{`U}
\DeclareMathSymbol{V}{\mathalpha}{operators}{`V}
\DeclareMathSymbol{W}{\mathalpha}{operators}{`W}
\DeclareMathSymbol{X}{\mathalpha}{operators}{`X}
\DeclareMathSymbol{Y}{\mathalpha}{operators}{`Y}
\DeclareMathSymbol{Z}{\mathalpha}{operators}{`Z}
\DeclareMathSymbol{e}{\mathalpha}{operators}{`e}

\newcommand{\mathbb}{\mathbf}

\newcommand{\RR}{\mathbb{R}}

\gdef\thechapter{\@Roman\c@chapter}%

\newtheorem{theoreme}{\noindent {\textbf{Th\'eor\`eme}}}[section]
{\nolinebreak  \end{trivlist}}

\renewcommand{\bar}{\overline}

\usepackage{epigraph}

\title{Le diamètre et la traversale~:~dans l'atelier de Girard Desargues}
\author{
par Marie Anglade et Jean-Yves Briend
}
\date{Decembre 12th, 2018}

\begin{document}
\maketitle

\setlength{\epigraphwidth}{11cm}
\epigraph{Le discours mathématique est une pensée qui veut se faire comprendre d'une autre pensée en lui devenant transparente.}{Vladimir Jankélévitch, \textit{La Musique et l'Ineffable}}

\begin{abstract} Le \textit{Brouillon Project} de Girard Desargues sur les coniques développe, dans sa partie centrale, la notion de \textit{traversale,} notion qui généralise celle de diamètre d'Apollonius et permet d'unifier le traitement des diverses espèces de coniques. Il est souvent écrit qu'il s'agit là d'un équivalent de la \textit{polaire,} concept émergeant au début du XIX${}^e$ siècle. Nous allons dans cet article explorer en détail les passages du texte de Desargues qui traitent de la traversale et de ses propriétés et montrer comment, en tenant compte des notes ajoutées postérieurement à la rédaction du premier jet du \textit{Brouillon,} on peut y reconstituer la naissance d'une nouvelle théorie, à savoir une \textit{théorie projective de la polarité associée à une conique.} Nous verrons qu'outre une correspondance naturelle entre points et droites, Desargues a compris qu'une conique induisait, sur chaque droite du plan, une involution et nous montrerons comment il vient progressivement à bout des difficultés conceptuelles considérables auxquelles il doit faire face pour exprimer clairement ses idées novatrices. 
\end{abstract}

\vspace*{.5cm}
\noindent{\bf Abstract.} In his \textit{Brouillon Project} on conic sections, Girard Desargues studies the notion of \textit{traversale}, which generalizes that of diameter introduced by Apollonius. One often reads that it is equivalent to the notion of \textit{polar,} a concept that emerged in the beginning of 19th century. In this article we shall study in great detail the developments around that notion in the middle part of the \textit{Brouillon project.} We shall in particular show, using the notes added by Desargues after the first draft was written, how here arises a new \textit{theory,} that of the \textit{polarity} associated to a conic section. We shall show that besides a natural correspondance between points and lines, Desargues has understood that a conic induces on every line in its plane an involution, and how he progressively handles the considerable conceptual difficulties he has to deal with when trying to explain clearly his novative ideas. 

\vspace*{1cm}

\section*{Introduction}
On retient couramment du \textit{Brouillon Project\footnote{Nous faisons ici référence au \textit{Brouillon project d'une atteinte aux evenemens des rencontres du Cone avec un Plan,} de Girard Desargues, paru en 1639, et nous basons sur l'unique exemplaire connu, conservé et numérisé par la Bibliothèque nationale de France, au département de la Réserve des livres rares, sous la référence RESM-V-276.}} sur les coniques de Girard Desargues son aspect novateur, tant dans la notion d'involution que dans les principes, que l'on peut qualifier de \textit{projectifs,} qui guident la démarche de l'auteur. On cite souvent comme résultat principal de ce texte le \textit{théorème d'involution} pour les pinceaux de coniques. Son énoncé et sa démonstration arrivent à la moitié environ du \textit{Brouillon,} après un long développement sur la notion d'involution. Il est légitime de se demander ce qui peut en constituer la seconde moitié et nous allons nous pencher sur cette question en examinant la notion arguésienne de \textit{traversale} qui permettra à Desargues d'unifier et de généraliser la théorie des diamètres d'Apollonius. Les commentateurs de ce texte\footnote{On pense en particulier à René Taton, dont l'édition dans \cite{taton} des {\oe}uvres mathématiques de Desargues fait référence, mais aussi au très joli article \cite{cracewski} de Stanislaus Chrzaszczewski, certes plus ancien.}, au courant de tout ce que la théorie des coniques doit à cette notion, reconnaissent dans la traversale ce que l'on nomme aujourd'hui la \textit{polaire.} Cependant leur analyse ne va guère plus loin qu'un simple constat et aucun n'a vraiment pris le soin de considérer le développement des idées, arguments et concepts dans le \textit{Brouillon Project.} Or ces développements permettent de comprendre qu'au delà de la notion de polaire, qui donne une correspondance entre points et droites, le \textit{Brouillon} contient en fait une théorie complète de la polarité, comme nous allons tenter de le montrer dans cet article.


Comme en d'autres instances dans le \textit{Brouillon Project,} le cheminement de Desargues est tortueux et parfois très difficile à suivre. Dans le présent article, nous chercherons à le clarifier. Desargues commence par définir la traversale d'un point relativement à une «~figure NB,NC~» qui n'est \emph{a priori} qu'une union de deux droites, dans une énonciation qui rappelle fortement celle de certaines propositions du livre VII de la \textit{Collection Mathématique} de Pappus. Il déclare alors immédiatement que la même notion est en fait attachée à \textit{une} conique donnée (passant par quatre points situés sur les deux droites NB,NC), mais il se rendra compte à la relecture de son texte, comme on peut s'en apercevoir en consultant les \textit{Advertissements,} notes d'ajouts et de corrections que l'on trouve annexées au \textit{Brouillon Project,} que cela n'est en rien une évidence.


Il revient, à la page 21, à la notion de traversale, dont il va démontrer qu'elle peut effectivement être rattachée à la donnée d'une conique et d'un point pris hors de celle-ci. On trouve dans ce passage des énoncés concernant des divisions harmoniques induites par des coniques qui sont très proches de ceux que l'on peut lire chez Apollonius, mais Desargues utilise pleinement la richesse de sa notion d'involution, tant dans son caractère agrégatif que dans ses relations subtiles avec la division harmonique\footnote{\textit{Voir} \cite{anglade-briend-1} pour plus de précision, ainsi que la section sur le vocabulaire arguésien.}, pour arriver finalement à ce qu'il convient d'appeler une théorie complète de la \textit{polarité.} 

Rappelons que cette théorie, outre qu'elle associe à une conique une correspondance entre points et droites (ici représentée par l'association, à tout point du plan, de sa traversale), montre également qu'une conique induit, dans toute droite du plan\footnote{Sauf celles qui sont tangentes à la conique.}, une correspondance involutive entre les points de cette droite (ici représentée par la notion d'involution de Desargues). Il y a dans ce passage une ambiguïté subtile quant à cette involution induite. En effet, pour construire la traversale, Desargues utilise son théorème d'involution. Or celui-ci fournit sur chaque droite une involution qui est différente de celle donnée par la conique, l'involution de polarité. Même s'il n'explicite pas vraiment le fait qu'il passe de l'une à l'autre, ce qui est nécessaire pour fonder sa théorie,  nous verrons que c'est bien cette dernière involution, induite par la conique, qu'il utilise. Desargues a en fait une vue très claire de ces subtilités, et nous pouvons envisager de le déclarer \textit{inventeur et pionnier de la théorie de la polarité.}


On peut se demander pourquoi cette théorie n'a pas survécu à Desargues et a dû attendre le début du XIX$^e$ siècle pour être redécouverte et exploitée avec le profit que l'on sait. 

Tout d'abord, le \textit{Brouillon Project} est un \textit{premier jet} de rédaction d'un traité, où Desargues commence à mettre de l'ordre dans ses idées, très novatrices, pour les soumettre à l'avis des géomètres de son temps. La lecture attentive du texte, en tenant compte des ajouts ultérieurs à sa première rédaction contenus dans les \textit{Advertissements,} peut laisser penser qu'une des raisons qui ont conduit Desargues à ne pas mener à bien une rédaction finale claire et compréhensible tient à la difficulté de la tâche~:~présenter une théorie comprenant une correspondance points-droites (une dualité) ainsi que l'induction par une conique dans chacune des droites du plan d'une involution, représente, au dix-septième siècle, un défi redoutable. Il suffit au lecteur pour s'en convaincre de lire l'appendice du présent article consacré à une présentation moderne de cette théorie~:~même avec l'aide de la théorie des ensembles, de l'algèbre linéaire et de la géométrie projective, la mise en ordre des arguments demande un grand soin et ne va pas de soi \textit{a priori.} 

De plus, les réactions peu enthousiastes (on peut penser à la tiédeur d'un Descartes, \textit{voir} \cite{taton}, p. 185) voire franchement hostiles (ici c'est Jean de Beaugrand et ses \textit{Advis Charitables} qui viennent à l'esprit, \textit{voir} \cite{anglade-briend-1}) qu'il a reçues l'ont sans doute découragé de mettre le tout  au propre dans une présentation qui permette d'emporter l'adhésion. 

Enfin, on ne peut négliger le fait que le développement rapide des méthodes analytiques issues entre autres des travaux de Descartes ou Fermat vont accaparer l'attention des mathématiciens pendant les décennies qui suivent la rédaction du \textit{Brouillon Project,} les détournant de recherches en «~géométrie pure~», considérée alors comme un domaine archaïque.

Le \textit{Brouillon Project} peut donc être comparé à une «~prépublication~» d'aujourd'hui. Nous pourrions définir ainsi le projet de cet article~:~donner une forme publiable à ce texte en en clarifiant l'exposition et en en comblant les lacunes. Nous tâchons pour cela de rester fidèles aux manières de faire de Desargues, classiquement euclidiennes. Nous commençons par quelques rappels sur le vocabulaire arguésien et sur la notion d'involution. Nous présentons ensuite brièvement le déroulement de l'original concernant la traversale et les diamètres, espèrant guider le lecteur dans un texte parfois déroutant. En nous basant sur une analyse détaillée des pages 21 et 22 du \textit{Brouillon} et en tenant compte de la chronologie du développement de ses idées, lisible au travers des ajouts et des corrections donnés dans les \textit{advertissements,} nous montrons enfin pourquoi l'on peut considérer que Desargues a bien développé une \emph{théorie de la traversale et de l'involution de polarité.} Nous en donnons une illustration en analysant son traitement très élégant des asymptotes aux hyperboles. Dans l'appendice de cet article, nous donnons une courte présentation, dans le langage contemporain de l'algèbre linéaire et de la géométrie projective, de la polarité associée à une conique. Cela nous a semblé utile car la théorie des coniques a depuis quelques décennies presqu'entièrement disparu des cursus de mathématiques, du moins en France. 

\vs
\noindent Nous voudrions faire de cet article un hommage au génie de Girard Desargues. 

\vs
\section*{Introduction: english version}
Of the various aspects of Desargues' \textit{Brouillon Project\footnote{We refer here to the \textit{Brouillon project d'une atteinte aux evenemens des rencontres du Cone avec un Plan,} by Girard Desargues, published in 1639. We base our study on the unique known original, which is kept at the Bibliothèque nationale de France, in the Réserve des livres rares, under the reference RESM-V-276. It can be downloaded from the Gallica web site.},} one often retains the notion of involution, its main theorem on pencils of conics, and the views the author adopts on geometry, which are very innovative and can be decribed as being \textit{projective.} But the heart of the \textit{Brouillon} deals with a notion introduced and extensively used by Desargues, that of \textit{traversale} of a point associated to a conic section, which allows to generalize and unify the apollonian notion of diameter. It is often written\footnote{One thinks for instance about René Taton, whose edition \cite{taton} of Desargues' mathematical works still is a reference, see also \cite{fieldgray}, or the very nice article \cite{cracewski} by Stanislaus Chrzaszczewski, although it is older.} that this notion is equivalent to what we call now the \textit{polar,} but no one has really gone further in analyzing the development of ideas, arguments and concepts in the text. We shall show, by a careful analysis of the original of the \textit{Brouillon,} taking into account the notes added afterwards by Girard Desargues, how he comes to realize that he has just built a complete theory of polarity.

The way the \textit{Brouillon Project} is written makes it sometimes hard to follow its author. We shall try to clarify this in the present article. The \textit{traversale} of a point is first defined, p. 10, with respect to a "figure" $NB,NC$, that is a union of two lines, in a way that strongly recalls some of the propositions of the 7th book of Pappus' \textit{Synagoge.} He immediately writes that it is in fact defined with respect to a conic, but he realizes that it is not at all obvious, as can be seen by reading his notes and addenda, entitled "Advertissements", written afterwards. So after having given this incomplete definition, Desargues states, proves and develops his theorem on involutions defined by a pencil of conics and he comes back to the traversale on p. 21. He states and proves various statements that strongly recall similar propositions regarding harmonic divisions defined by conics thta can be found in Apollonius. But he strongly uses the full force of his notion of involution, including its "agregative" character\footnote{see \cite{anglade-briend-1} for a thorough investigation of Desargues' notion of involution.} to finally build a complete theory of polarity. 

Let us recall that this theory defines a correspondance between lines and points, but also shows that on each line of a plane, a conic induces an involution. One could recover such a statement from the propositions of Apollonius on harmonic divisions induced by a conic on a line, but one has to realize that for this the line \textit{has to cut the conic.} Desargues shows, using his stronger notion of involution, that a conic induces an involution on every\footnote{Except the lines that are tangent to the conic.} line of the plane, even those which do not cut the conic. There is nevertheless an ambiguity on this involution. To define the traversale, Desargues uses his theorem on involution induced by a \textit{pencil} of conics through for base points. The involution thus defined is (almost) never equal to the one defined by polarity of a single conic in the pencil, and Desargues switches from the first to the second in a very discreet way, but, as the rest of the text shows, he has a very clear view on those matters. Following Poncelet, we can thus christen Desargues the \textit{father and founder of the projective theory of polarity.}

One may wonder why this wonderful theory has been forgotten almost immediately after the publication of the \textit{Brouillon.} 

Firstly, it is quite obvious that the \textit{Brouillon Project} is really a first draft, wher Desargues tries to expose as clearly as he can the present early stage of his thoughts, to submit them to the geometers of his time. Our careful reading of the \textit{Brouillon,} including the "dynamical" view that one gets by considering the addenda contained in the \textit{Advertissements} shows that Desargues has faced considerable conceptual difficulties~:~prensenting a complete theory of polarity in the first half of 16th century is a challenge that would have turned away the most capable of mathematicians. Even with the full force of moderne projective theory, including concepts of set theory and linear algebra, a clear and concise presentation of this topic does not come without some substantial effort. 

Secondly, his contemporaries had quite mitigated reacions to the \textit{Brouillon,} and some of them even showed some clear hositility (see for instance Jean de Beaugrand in his "Advis Charitables", as commented upon in \cite{anglade-briend-1}). 

And finally we can not neglect the fact that "pure" geometry became rapidly unfashionnable in the mathematical community as the development of analytic methods, following the publication of Descartes' \textit{Geometry} in 1637.

To conclude, we can compare the \textit{Brouillon Project} with today's preprints. Our hope in this article is to make its content far more accessible to the reader than it is in the original. We shall begin by short reminders about Desargues' vocabulary and his notion of involution, then we will give a short description of the way the arguments are arranged in the original and immediately afterwards we shall give a precise analysis of those developments, completing and sorting things to show how Desargues really builds a new theory, that of a polarity induced by a conic section. Before our conclusion, we will show how Desargues applies his theory of the traversale in a very elegant manner to treat the asymptots of an hyperbola. In a appendix, we give a short description of polarity theory in the modern languages of projective geometry, as it is not certain that every reader has acquaintances with that topic.

We would like to see this article as a tribute to the genius of Girard Desargues.


\section{Le vocabulaire arguésien et la notion d'involution}
L'un des principaux reproches faits à Desargues par ses contemporains est l'emploi d'une terminologie qui lui est propre et qui déroge aux usages de son temps. De celle-ci il ne nous reste aujourd'hui que le terme d'\textit{involution} et il a fallu attendre le vingtième siècle pour que des mathématiciens, ceux du groupe Bourbaki par exemple, se remettent à déployer un vocabulaire imagé, d'inspiration sensible, dans le champ mathématique\footnote{À propos de l'aspect \textit{baroque} du vocabulaire arguésien, voir l'article \cite{mesnard} de Jean Mesnard.}.

Le \textit{Brouillon Project} commence par une mise au point sur ce que Desargues entend par ligne et plan, qu'il conçoit comme étendus «~au besoin à l'infini~» de toutes parts\footnote{p. 1, l. 11 et l. 30}. Il définit ensuite la notion d'\textit{ordonnance,} comme étant la donnée d'un certain nombre de droites toutes concourantes en un point donné ou toutes parallèles entre elles. Il décide d'unifier les deux notions de concourance et de parallélisme en introduisant un point de concours, qu'il nomme \textit{but de l'ordonnance} et qui se situe à distance finie ou infinie selon les cas\footnote{p. 1, l. 15 à 29}. Il fait de même pour une famille de plans\footnote{p. 1, l. 33 à 46}, la droite de concours étant ici nommée \textit{essieu\footnote{Notons que ce terme est celui utilisé en gnomonique.}} de l'ordonnance, en introduisant au passage ce qu'il décrira plus loin comme une ligne dont tous les points sont à distance infinie. 

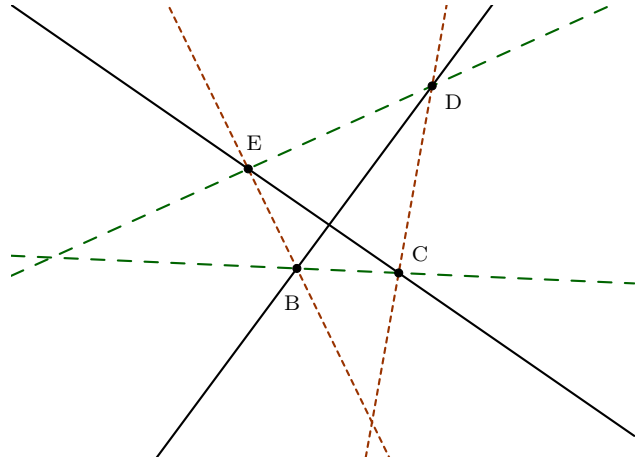
\begin{figure}[!ht]
\centering
\definecolor{zzttqq}{rgb}{0.6,0.2,0.}
\definecolor{qqwuqq}{rgb}{0.,0.39215686274509803,0.}
\begin{tikzpicture}[line cap=round,line join=round,>=triangle 45,x=1.0cm,y=1.0cm]
\clip(-3.66,-2.08) rectangle (4.54,3.9);
\draw [line width=0.8pt,domain=-3.66:4.54] plot(\x,{(--2.7-1.38*\x)/1.98});
\draw [line width=0.8pt,domain=-3.66:4.54] plot(\x,{(--0.5056--2.42*\x)/1.78});
\draw [line width=0.8pt,dash pattern=on 5pt off 5pt,color=qqwuqq,domain=-3.66:4.54] plot(\x,{(-4.8048-1.1*\x)/-2.42});
\draw [line width=0.8pt,dash pattern=on 5pt off 5pt,color=qqwuqq,domain=-3.66:4.54] plot(\x,{(-0.5688--0.06*\x)/-1.34});
\draw [line width=0.8pt,dash pattern=on 2pt off 2pt,color=zzttqq,domain=-3.66:4.54] plot(\x,{(--3.4128-2.48*\x)/-0.44});
\draw [line width=0.8pt,dash pattern=on 2pt off 2pt,color=zzttqq,domain=-3.66:4.54] plot(\x,{(--0.4008-1.32*\x)/0.64});
\begin{scriptsize}
\draw [fill=black] (-0.54,1.74) circle (1.5pt);
\draw[color=black] (-0.46,2.09) node {$E$};
\draw [fill=black] (0.1,0.42) circle (1.5pt);
\draw[color=black] (0.04,-0.03) node {$B$};
\draw [fill=black] (1.44,0.36) circle (1.5pt);
\draw[color=black] (1.72,0.61) node {$C$};
\draw [fill=black] (1.88,2.84) circle (1.5pt);
\draw[color=black] (2.16,2.63) node {$D$};
\end{scriptsize}
\end{tikzpicture}
\caption{Un quadrangle de bornes $B,C,D,E$ et ses trois couples de bornales $BE,DC;BC,DE;BD,CE$.}\label{Bornales}
\end{figure}

Il donne un nom particulier aux sommets d'un quadrilatère en les appelant des \textit{bornes,} et aux six diagonales de ce quadrilatère qu'il nomme \textit{bornales} et que l'on peut naturellement regrouper en trois \emph{couples} de droites\footnote{\textit{p. 2, l. 51 à 56}} (\textit{voir} la figure \ref{Bornales}).

\begin{figure}[!ht]
\centering
\definecolor{qqqqcc}{rgb}{0.,0.,0.8}
\definecolor{uuuuuu}{rgb}{0.26666666666666666,0.26666666666666666,0.26666666666666666}
\definecolor{ffqqtt}{rgb}{1.,0.,0.2}
\definecolor{xdxdff}{rgb}{0.49019607843137253,0.49019607843137253,1.}
\definecolor{qqqqff}{rgb}{0.,0.,1.}
\begin{tikzpicture}[line cap=round,line join=round,>=triangle 45,x=0.73126142595978cm,y=0.73126142595978cm]
\clip(-4.3,-2.46) rectangle (6.64,6.3);
\draw [line width=1.6pt,domain=-4.3:6.64] plot(\x,{(--2.4188--2.08*\x)/2.82});
\draw [line width=1.6pt,color=ffqqtt] (-2.6,-1.06)-- (-0.36,5.02);
\draw [line width=2.pt,color=ffqqtt] (-2.6,-1.06)-- (0.22,1.02);
\draw [line width=1.6pt,color=ffqqtt] (2.3289004619907607,2.5755010499790005)-- (-3.32,4.26);
\draw [line width=1.6pt,color=qqqqcc] (-3.32,4.26)-- (-0.9052874475411561,3.5399340709597187);
\draw [line width=1.6pt,color=ffqqtt] (2.3289004619907607,2.5755010499790005)-- (5.386131037379254,4.8304796304073925);
\begin{scriptsize}
\draw [fill=qqqqff] (-2.6,-1.06) circle (2.5pt);
\draw[color=qqqqff] (-3.06,-0.79) node {$H$};
\draw [fill=qqqqff] (0.22,1.02) circle (2.5pt);
\draw[color=qqqqff] (-0.06,1.35) node {$B$};
\draw[color=black] (-3.82,-1.35) node {tronc};
\draw [fill=xdxdff] (2.3289004619907607,2.5755010499790005) circle (2.5pt);
\draw[color=xdxdff] (2.92,2.43) node {nœud};
\draw [fill=xdxdff] (5.386131037379254,4.8304796304073925) circle (2.5pt);
\draw[color=xdxdff] (5.04,5.15) node {$G$};
\draw [fill=qqqqff] (-0.36,5.02) circle (2.5pt);
\draw[color=qqqqff] (-0.22,5.39) node {$A$};
\draw[color=ffqqtt] (-1.94,1.99) node {rameau déployé};
\draw[color=ffqqtt] (0.44,-0.19) node {rameau plié au tronc};
\draw [fill=qqqqff] (-3.32,4.26) circle (2.5pt);
\draw[color=qqqqff] (-3.18,4.63) node {$D$};
\draw [fill=uuuuuu] (-0.9052874475411561,3.5399340709597187) circle (1.5pt);
\draw[color=uuuuuu] (-0.64,3.79) node {$d$};
\draw[color=qqqqcc] (-2.98,3.63) node {brin de rameau};
\end{scriptsize}
\end{tikzpicture}
\caption{Un tronc, des n{\oe}uds, des rameaux et des brins de rameaux, pliés ou déployés.}\label{Rameaux}
\end{figure}
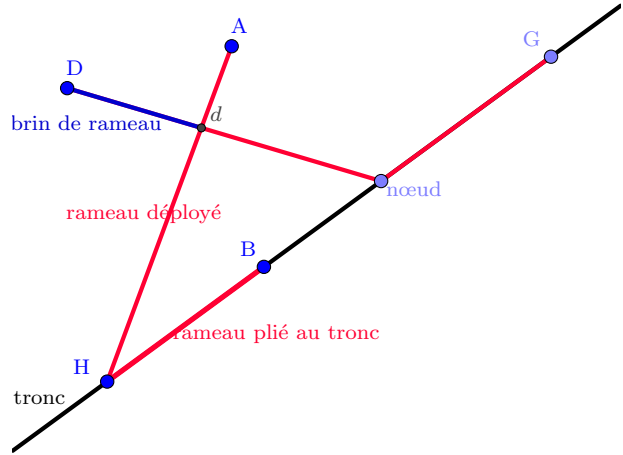

Le texte se poursuit, jusqu'à la page 10, par un long développement sur la notion d'involution, configuration particulière de trois couples de points situés sur une même droite. Lorsque Desargues s'intéresse à des points alignés sur une droite qu'il veut mettre en valeur, il parle de cette droite comme d'un \textit{tronc.} Les points en question sur cette droite sont qualifiés de \textit{n{\oe}uds.} De ces n{\oe}uds peuvent partir des segments de droites que l'auteur nomme des \textit{rameaux.} Si ces rameaux sont sur le tronc, il sont \textit{pliés au tronc,} et si au contraire ils lui sont transverses, alors ils sont \textit{déployés au tronc.} Deux tels rameaux peuvent se couper, et les nouveaux segments ainsi définis dans les rameaux s'appellent des \textit{brins de rameaux\footnote{\textit{voir} la p.2 de l'original.},} comme il est illustré sur la figure \ref{Rameaux}.

Desargues définit ensuite certaines conditions d'ordonnancement de couples de n{\oe}uds sur un tronc, en disant par exemple que $A$ est \emph{engagé} entre les deux segments $AB,AC$ si $A$ est situé entre $B$ et $C$, et qu'il est \emph{dégagé} sinon. Ou bien encore que deux segments (ou couples de points) $BH,CG$ sont \emph{mêlés} si une et une seule des extrémités de l'un des deux segments est située entre les deux extrémités de l'autre, et qu'ils sont \emph{démêlés} sinon. Il définit alors un \textit{arbre} comme la donnée d'un point $A$ qu'il nomme souche et de trois couples de n{\oe}uds $B,H;C,G;D,F$ de sorte que d'une part $A$ soit dans la même situation d'engagement ou de dégagement vis-à-vis des couples de branches $AB,AH;AC,AG;AD,AF$ et que d'autre part les égalités de rectangles suivantes aient lieu~:~
\[
AB.AH=AC.AG=AD.AF.
\] 
On dit par ailleurs que les trois couples $B,H;C,G;D,F$ forment une \textit{involution\footnote{Ou~«~sont en involution~»}} si d'une part ils sont tous de même \textit{mêlés} ou \textit{démêlés} les uns aux autres et si d'autre part les égalités de rapports suivantes sont vérifiées~:~
\[
\frac{GD.GF}{CD.CF}=\frac{GB.GH}{CB.CH},\; \frac{FC.FG}{DC.DG}=\frac{FB.FH}{DB.DH},\;\mbox{et}\;\frac{HC.HG}{BC.BG}=\frac{HD.HF}{BD.BF}.
\]
Les rapports ci-dessus sont construits à partir d'une combinatoire des rectangles (les \emph{relatifs} et les \emph{gémeaux}) mise au point par Desargues et nous renvoyons à la figure \ref{involution}, qui en donne une illustration pour le premier rapport de la série. Desargues démontre\footnote{Pour plus de détails sur tout ce qui concerne les arbres et les involutions, nous renvoyons à l'article \cite{anglade-briend-1} des auteurs du présent texte.} alors que la configuration $A;B,H;C,G;D,F$ est un arbre si et seulement si la configuration $B,H;C,G;D,F$, d'où la souche a disparu, forme un arrangement de trois couples de n{\oe}uds en involution.

\begin{figure}[!ht]
\centering
\begin{tikzpicture}
\draw (1,0) node[above] {$GF.GD$};
\draw (0,0) -- (2,0);
\draw (1,0) node[below] {$CF.CD$};
\draw [<->] (0,-0.3) to[bend left] (0,0.3);
\draw (-1,0) node {relatifs};
\draw (2.2,0) node {$=$};
\draw (2.5,0) -- (4.5,0);
\draw (3.5,0) node[above]{$GB.GH$};
\draw (3.5,0) node[below]{$CB.CH$};
\draw [<->] (4.5,-0.3) to[bend right] (4.5, 0.3);
\draw (5.4,0) node {relatifs};
\draw  [<->] (1,0.7) to[bend left] (3.5,0.7);
\draw (2.3, 1.4) node {gémeaux};
\end{tikzpicture}
\caption{Combinatoire de rectangles~:~rapports de rectangles \emph{relatifs} égaux aux rapports de leurs \emph{gémeaux}.}\label{involution}
\end{figure}

Il n'est pas inutile de rappeler ici comment on peut interpréter dans le langage contemporain la notion d'involution~:~les trois couples de points $B,H;C,G;D,F$, situés sur une droite projective $\Delta$ sont en involution si et seulement s'il existe une homographie involutive\footnote{Elle est égale à sa réciproque.} de $\Delta$ sur elle-même qui envoie chaque n{\oe}ud sur son accouplé\footnote{\textit{Voir} par exemple le livre \cite{sidler} pour de plus amples détails sur les involutions d'une droite projective.}.

Dans son étude de la notion d'involution, Desargues passe un temps certain à analyser le cas particulier où des points viennent à se confondre. Si dans un arbre $A;B,H;C,G;D,F$ un rectangle comme $AC,AG$ est en fait un carré, il déclare que les n{\oe}uds  $C,G$ sont des \textit{n{\oe}uds moyens.} Deux cas se présentent alors. 

Si $C=G$, Desargues dit que $C,G$ forment un couple de n{\oe}uds moyens \emph{doubles} et il  démontre que l'involution en admet nécessairement un autre comme $D,F$ par exemple, de sorte que la souche $A$ de l'arbre associé est milieu du segment $CD$. Desargues dit alors que $B,H;C,C;D,D$ est une involution \emph{de quatre points seulement,} et $B,H,C,G$ forment en fait une \emph{division harmonique,} notion qui n'est jamais nommée ni mentionnée en tant que telle dans le \textit{Brouillon Project.} Les n{\oe}uds moyens doubles correspondent, dans l'homographie involutive définie par l'involution de Desargues, à deux points fixes. Cette homographie est donc hyperbolique dans ce cas.


Si au contraire $C\neq G$, de sorte que $A$ est le milieu de $CG$, les deux points $C$ et $G$ sont dits être des n{\oe}uds moyens \textit{simples.} L'involution a alors deux couples de n{\oe}uds moyens simples, à savoir $C,G$ et $G,C$. Dans ce cas, l'homographie involutive correspondante est sans point fixe, ou elliptique, et la définition de n{\oe}ud moyen n'a de sens qu'eu égard au choix de la souche $A$; il s'agit donc d'une notion \emph{affine\footnote{\textit{Voir} l'article \cite{anglade-briend-1}.}.} 

Dans une situation où l'on a identifié des n{\oe}uds moyens, les n{\oe}uds de l'involution qui ne le sont pas sont dits \textit{extrêmes.} 

La notion de n{\oe}ud moyen double joue un rôle éminent dans le \textit{Brouillon Project,} comme nous le verrons dans la suite de cet article. Rappelons que si $B,H;C,G$ est une division harmonique, on peut en déduire \emph{deux} involutions (\textit{voir} par exemple \cite{sidler}, p. 27), en précisant quel couple de n{\oe}uds on choisit pour être les moyens doubles~:~
\[
B,B;H,H;C,G \;\mbox{ou}\;B,H;C,C;G,G.
\]
Les  deux arbres correpondant à ces involutions ont des souches dites \emph{réciproques.} Ainsi, la notion d'involution généralise et précise celle de division harmonique.  

Mentionnons pour finir un cas particulier de la division harmonique ou de l'involution à n{\oe}uds moyens doubles que Desargues traite précisément et qui aura son importance par la suite. Si quatre points $B,H;C,C;G,G$ sont en involution, alors $B$ est le milieu du segment $CG$ si et seulement si $H$ est le point à l'infini sur la droite.

\begin{figure}[!ht]
\centering
\definecolor{wwccff}{rgb}{0.4,0.8,1.}
\begin{tikzpicture}[line cap=round,line join=round,>=triangle 45,x=0.7317073170731712cm,y=0.7317073170731712cm]
\clip(-3.14,-8.02) rectangle (13.26,0.94);
\draw [line width=1.6pt,domain=-3.14:13.26] plot(\x,{(-88.8288-0.04*\x)/14.});
\draw [domain=-3.14:13.26] plot(\x,{(-45.9248--3.82*\x)/8.28});
\draw [domain=-3.14:13.26] plot(\x,{(-39.01343157894736--3.827894736842105*\x)/5.516842105263157});
\draw [domain=-3.14:13.26] plot(\x,{(-33.565638651854144--3.834117656661959*\x)/3.3388201683141636});
\draw [domain=-3.14:13.26] plot(\x,{(-30.105445614035087--3.8380701754385966*\x)/1.955438596491228});
\draw [domain=-3.14:13.26] plot(\x,{(-15.719560233918116--3.8545029239766078*\x)/-3.7960233918128674});
\draw [domain=-3.14:13.26] plot(\x,{(-10.9072--3.86*\x)/-5.72});
\draw [line width=1.6pt,color=wwccff,domain=-3.14:13.26] plot(\x,{(--38.70211365408517-5.228587089942698*\x)/2.1072139708549855});
\begin{scriptsize}
\draw [fill=black] (-1.72,-6.34) circle (1.5pt);
\draw[color=black] (-1.86,-6.05) node {$F$};
\draw [fill=black] (12.28,-6.38) circle (1.5pt);
\draw[color=black] (12.42,-6.09) node {$H$};
\draw[color=black] (-2.68,-5.99) node {$\Delta$};
\draw [fill=black] (1.0431578947368425,-6.347894736842105) circle (1.5pt);
\draw[color=black] (0.82,-6.03) node {$D$};
\draw [fill=black] (3.221179831685836,-6.354117656661959) circle (1.5pt);
\draw[color=black] (3.02,-6.07) node {$B$};
\draw [fill=black] (4.604561403508772,-6.358070175438597) circle (1.5pt);
\draw[color=black] (4.34,-6.07) node {$C$};
\draw [fill=black] (10.356023391812867,-6.374502923976608) circle (1.5pt);
\draw[color=black] (10.44,-6.15) node {$G$};
\draw [fill=black] (6.56,-2.52) circle (1.5pt);
\draw[color=black] (5.94,-2.49) node {$K$};
\draw [fill=wwccff] (7.597183243619282,-0.48425104785043516) circle (1.5pt);
\draw[color=black] (7.60,-0.03) node {$c$};
\draw [fill=wwccff] (9.704397214474268,-5.712838137793133) circle (1.5pt);
\draw[color=black] (9.36,-5.79) node {$g$};
\draw[color=black] (7.35,0.75) node {$\delta$};
\draw [fill=wwccff] (7.829908323381838,-1.061707046331059) circle (1.5pt);
\draw[color=black] (7.88,-0.67) node {$b$};
\draw [fill=wwccff] (8.01168299928028,-1.5127408458533211) circle (1.5pt);
\draw[color=black] (8.08,-1.19) node {$d$};
\draw [fill=wwccff] (8.12638348273628,-1.7973448183511358) circle (1.5pt);
\draw[color=black] (8.41,-1.95) node {$f$};
\draw [fill=wwccff] (9.11156757621714,-4.241862035698981) circle (1.5pt);
\draw[color=black] (9.26,-3.95) node {$h$};
\end{scriptsize}
\end{tikzpicture}
\caption{La configuration de la \emph{ramée}}\label{Ramee-1}
\end{figure}
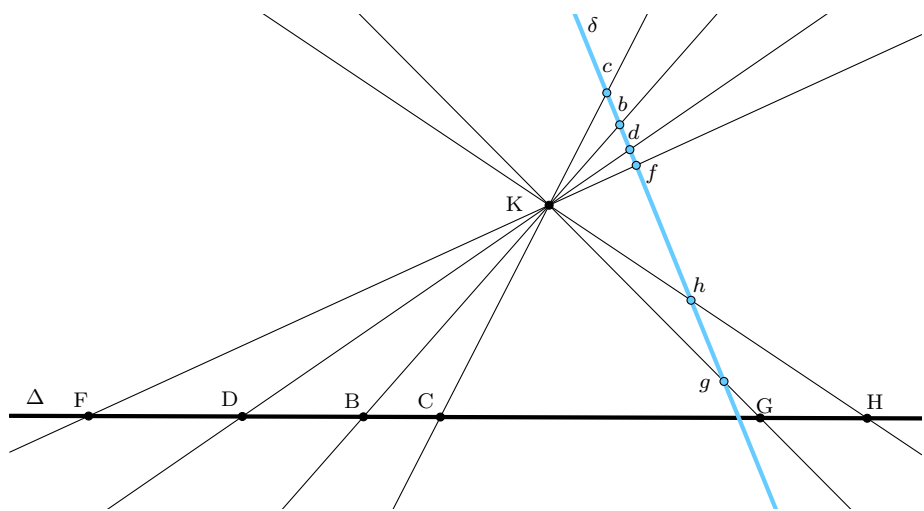

Ce n'est qu'à la page 11 que Desargues passe de la géométrie de la droite à celle du plan. Il y définit\footnote{p. 11, l. 9 à 14.} la notion de \textit{ramée} d'un arbre ou d'une involution, \textit{voir} la figure \ref{Ramee-1}. Si $B,H;C,G;D,F$ forment sur une droite $\Delta$ une involution, et si $K$ est un point pris hors de $\Delta$, alors les rameaux déployés $KH,KB;KC,KG;KD,KF$ forment la \textit{ramée} de l'involution $B,H$ \textit{etc.} Il démontre derechef le \textit{théorème de la ramée}~:~si $\delta$ est une droite du plan de la ramée, les rameaux de ladite ramée coupent $\delta$ en des points $b,h;c,g;d,f$ qui, eux aussi, forment une involution\footnote{\textit{Voir} l'article \cite{anglade-briend-2} pour plus de détails sur ce théorème et sa démonstration.}. Les pages 11 à 13 sont dévolues à la démonstration de ce théorème et à l'étude de certains de ses cas particuliers (parallélisme de certains rameaux, cas d'involution de quatre points).

À partir de la page 14 du \textit{Brouillon Project,} Desargues s'attaque à ce qui en forme le sujet principal~:~les coniques. Fidèle à son tropisme baroque, il introduit une terminologie entièrement nouvelle, ce qui provoquera quelques sarcasmes chez ses contemporains comme Jean de Beaugrand, pour qui ces termes~«~sont plus capables de mettre les esprits {\it en involution} ou d'en faire des {\it souches réciproques,} que de leur donner quelque nouvelle lumière dans les Mathématiques\footnote{{\it Voir} le fascicule côté V-12591 de la Bibliothèque nationale de France contenant, après les \textit{Advis Charitables} attaquant les travaux de Desargues, une longue lettre de Beaugrand critiquant le \textit{Brouillon Project.}}~». 

Étant donnés un cercle et un point $S$ hors du plan du cercle, la surface engendrée par le mouvement de droites issues de $S$ et passant par un point variable sur le cercle est nommée \textit{surface de rouleau\footnote{p. 14, l. 1 à 5.}.} Le rouleau proprement dit est le volume enserré par cette surface (notons que c'est ce volume que Desargues définit en premier). Le point $S$ en est le \textit{sommet\footnote{p. 14, l. 6.}} et le cercle servant à le former en est appelé la \textit{plate assiette\footnote{p. 14, l. 7. Il s'agit ici encore d'une terminologie issue de la gnomonique.}} ou \textit{base.} Si le sommet est à distance infinie, le rouleau est un \textit{cylindre\footnote{p. 14, l. 13.},} tandis que s'il est à distance finie on aura affaire à un \textit{cornet\footnote{p. 14, l. 15.}.} Nous emploierons dans la suite de cet article une terminologie plus en rapport avec l'usage, en nommant \emph{c\^one} la surface réglée engendrée par un cercle et un point, que ce point soit à distance finie ou infinie, en précisant le cas échéant qu'il peut s'agir d'un cylindre.

\begin{figure}[!ht]
\centering
\includegraphics[width=7cm]{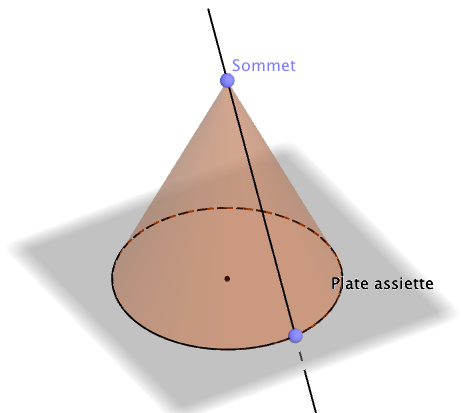}
\caption{Un rouleau, son sommet et sa plate assiette.}\label{Cone-01}
\end{figure}

Un plan différent de celui de la plate assiette et qui vient couper un rouleau est simplement nommé \textit{plan de coupe\footnote{p. 14, l. 25.}} par Desargues. La portion de surface de ce plan découpée par le rouleau s'appelle une \textit{coupe de rouleau\footnote{p. 15, l. 17.}} et la courbe qui forme son bord, une conique, est nommée par Desargues \textit{bord\footnote{p. 15, l. 20.}} de la coupe de rouleau. Il donne alors aux trois cas usuels de l'ellipse, de la parabole et de l'hyperbole les noms respectifs de \textit{défaillement,} d'\textit{égalation} et  d'\textit{outrepassement\footnote{p. 15, l. 24 à 32.}.} Ce n'est finalement qu'une tentative de traduction en français courant de la terminologie introduite par Apollonius, qui répond au souci pédagogique de Desargues de toucher un public plutôt concerné par les applications de la géométrie (perspective, coupe des pierres \textit{etc.}) que par cette discipline en tant que telle\footnote{\textit{Voir} à ce sujet la lettre de Descartes à Desargues, p. 185 de \cite{taton}.}. Il est à noter que ces termes ne seront pas utilisés par Desargues dans la suite du texte et qu'il restera de ce point de vue fidèle à l'usage classique. Notons que le cas des coniques dégénérées, unions de deux droites, sont comprises par Desargues comme étant des coupes de rouleau.

Au delà des noms choisis, c'est la manière dont Desargues décrit chaque type de conique qui est significative et montre son aspiration à l'unification. La simple lecture de ces quelques lignes\footnote{p. 15, l. 22 à 32} extraites du \textit{Brouillon Project} en donnera une idée plus claire qu'un long discours~:~

--~«~Quand le bord d'une coupe de rouleau se trouve estre deux droictes, le but de leur ordonnance est à distance ou finie, ou infinie.

Quand le bord d'une coupe de rouleau se trouve estre une ligne courbe, laquelle à distance finie rentre \& repasse en soy mesme, la figure en est nommée ou {\it Cercle,} ou {\it Ovale,} autrement {\it Elipse,} en francez, deffaillement.

Quand le bord d'une coupe de rouleau se trouve estre une ligne courbe, laquelle à distance infinie rentre \& repasse en soy mesme, la figure en est nommée {\it Parabole,} en francez, égalation.

Quand le bord d'une coupe de rouleau se trouve estre une ligne courbe, laquelle à distance infinie se mypartit en deux moitiez opposez dos à dos, la figure en est nommée {\it Hyperbole,} en francez, outrepassement ou excedement.~»

On note en particulier que les deux branches d'une hyperbole sont conçues comme constituant un seul et même objet. 

Desargues achève ce passage p. 15 avec les définitions de \textit{traversale} et d'\textit{ordonnées} d'une traversale. L'un des buts de cet article est de clarifier ce qu'il entend par là et nous allons maintenant nous consacrer à l'analyse précise de ces notions, en commençant par décrire rapidement le déroulement de leur étude dans le \textit{Brouillon Project.}
\section{Déroulement des arguments dans l'original du \textit{Brouillon Project}}
Comme son titre l'indique et comme les auteurs du présent texte l'ont déjà souligné ailleurs, l'organisation du \textit{Brouillon Project} ne permet pas une progression linéaire dans les arguments mais va et vient, se répète, laisse parfois dans l'obscurité une transition essentielle, dans ce qu'il faut bien qualifier de labyrinthe baroque. En donnant une description rapide, sans rentrer dans les détails, de l'enchaînement des arguments dans le texte original, nous voudrions fournir au lecteur potentiel du \textit{Brouillon} comme un fil d'Ariane.


Rappelons que les premières pages sont dédiées à la mise en place d'une terminologie qui va permettre de développer un discours de nature projective concernant des configurations de points sur une droite ou de droites dans un plan, et que jusqu'à la page 10 de l'original Desargues développe la théorie de l'involution. À la fin de la page 10 il rappelle l'énoncé du théorème de Ménélaüs, qu'il démontre au début de la page 11. Il passe alors à la notion de ramée puis à son théorème dit de \textit{la ramée,} énonçant que la notion d'involution est \textit{invariante par perspective.} Il en explore divers cas particuliers jusqu'à la fin de la page 13 où il commence à introduire sa terminologie concernant les cônes et les coniques.

\begin{figure}[!ht]
\centering
\definecolor{uuuuuu}{rgb}{0.26666666666666666,0.26666666666666666,0.26666666666666666}
\begin{tikzpicture}[line cap=round,line join=round,>=triangle 45,x=1.0cm,y=1.0cm]
\clip(0.76,-4.44) rectangle (9.44,3.1);
\draw [rotate around={29.05460409907715:(5.25,-0.41)},line width=2.pt] (5.25,-0.41) ellipse (3.27475669150986cm and 1.7310781000834738cm);
\draw [line width=2.pt,domain=0.76:9.44] plot(\x,{(--2365.1696597410687-382.69780153459897*\x)/-868.3078089864007});
\draw [line width=2.pt,domain=0.76:9.44] plot(\x,{(--42.737015331259315-5.696594353538518*\x)/1.0854517956318521});
\draw [line width=2.pt,domain=0.76:9.44] plot(\x,{(--15.734995889604445-5.696594353538518*\x)/1.0854517956318521});
\begin{scriptsize}
\draw[color=black] (1.32,-1.85) node {$\delta$};
\draw[color=black] (6.86,2.95) node {$\omega$};
\draw [fill=uuuuuu] (7.168334965260984,1.752190978479834) circle (2.0pt);
\draw[color=uuuuuu] (7.38,2.07) node {$C$};
\draw [fill=uuuuuu] (7.631240253420358,-0.677197098783236) circle (2.0pt);
\draw[color=uuuuuu] (7.26,-0.39) node {$G$};
\draw [fill=uuuuuu] (7.39978760934067,0.537496939848299) circle (2.0pt);
\draw[color=uuuuuu] (7.64,0.95) node {$B$};
\draw [fill=uuuuuu] (2.8046058570279646,-0.22267778303449617) circle (2.0pt);
\draw[color=uuuuuu] (3.12,-0.19) node {$C'$};
\draw [fill=uuuuuu] (3.2493705677652858,-2.556861714587676) circle (2.0pt);
\draw[color=uuuuuu] (2.8,-2.77) node {$G'$};
\draw [fill=uuuuuu] (3.0269882123966245,-1.389769748811086) circle (2.0pt);
\draw[color=uuuuuu] (3.34,-1.57) node {$B'$};
\end{scriptsize}
\end{tikzpicture}
\caption{Une conique, un diamètre et deux de ses ordonnées}\label{Diametre-01}
\end{figure}
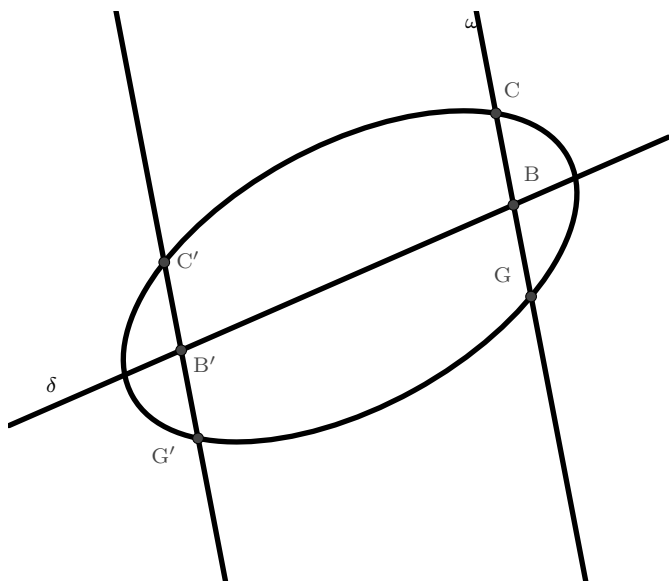

C'est à partir de la ligne 37 de la page 15 qu'il commence sa véritable exploration des rencontres d'un cône et d'un plan, en introduisant la notion de \textit{traversale.} Desargues semble vouloir généraliser la notion apollonienne de diamètre\footnote{Nous renvoyons à l'article \cite{hogendijk} de Jan Hogendijk pour une analyse des rapports en les conceptions apolloniennes et arguésiennes des coniques.}. Si $\cC$ est une section conique et si l'on se donne une \emph{direction,} c'est-à-dire, dans le langage arguésien, une ordonnance de droites parallèles entre elles, qu'Apollonius nomme dans ce contexte des \emph{ordonnées,} alors on peut trouver une droite $\delta$, dite \emph{diamètre,} associée à la famille des ordonnées, de sorte que les droites de l'ordonnance sont «~myparties~» par la conique et le diamètre, \textit{voir} la figure \ref{Diametre-01}. Dit autrement, si $\omega$ est une ordonnée quelconque qui coupe $\cC$ en les points $C$ \& $G$ et $\delta$ en $B$, alors $B$ est le milieu de $CG$. Les ordonnées sont «~ordonnées à un même but~» qui est le point à l'infini en chacune d'elles. Desargues va montrer que l'on peut ramener ce point dans le plan, à distance finie, et que les diamètres apolloniens ne sont qu'un cas particulier d'une notion plus générale, celle de \emph{traversale.}

La première définition qu'il en donne\footnote{p. 15, l. 37 à 42.} est quelque peu ambigüe, car elle se réfère à une «~figure $NB, NC$~», qui n'est \textit{a priori} qu'un couple de droites se coupant au point $N$, et nullement une conique générale, comme le laisse entendre Taton dans sa note de la page 138 de \cite{taton}. Vue son importance, citons-là \textit{in extenso~}:~ 

«~Quand en un plan une figure NB, NC, est rencontrée de plusieurs droites F C B, F I K, F X Y, d'une mesme ordonnance entre elles, \& qu'une mesme droicte N G H O, donne en chacune de ces droictes d'une mesme ordonnance entre elles un poinct G, H, O, couplé au but F, de leur ordonnance en involution avec les deux poincts comme X Y, I K, C B, qu'y donne le bord de la figure NB, NC, une telle droicte N G H, est pour cela nommée icy \textit{Traversale,} des droictes de l'ordonnance au but F, à l'égard de cette figure NB, NC, \& les droictes de cette ordonnance au but F, sont pour cela nommées \textit{Ordonnées} de la traversale N, G, H, à l'égard de la mesme figure NB, NC.~»

\begin{figure}[!ht]
\centering
\includegraphics[width=10cm]{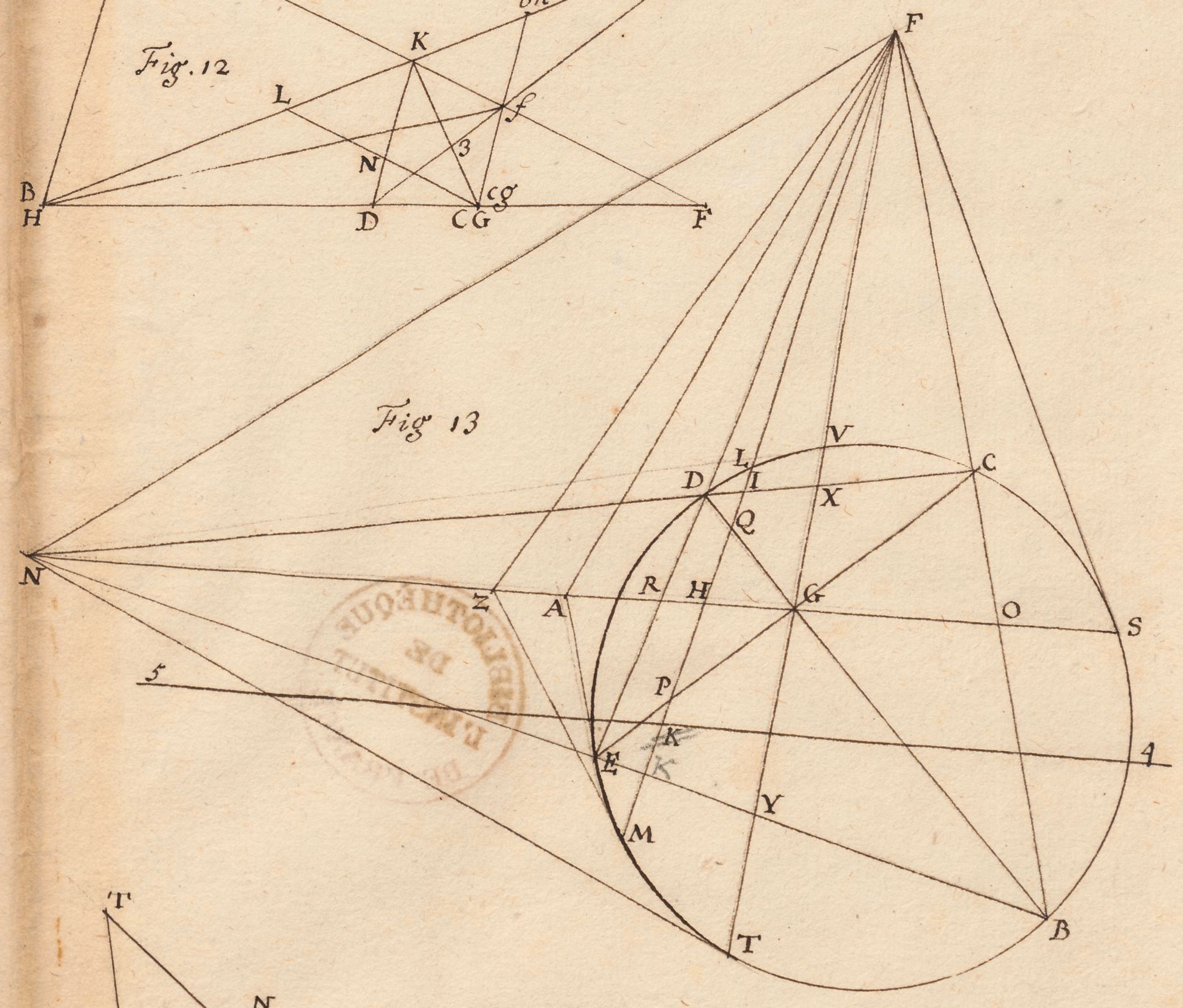}
\caption{La figure 13 du manuscrit de de la Hire illustrant la définition de la traversale.}\label{LaHire-01}
\end{figure}

La figure \ref{LaHire-01} montre le schéma qu'en donne Philippe de la Hire dans sa copie manuscrite du \textit{Brouillon Project\footnote{Ce manuscrit est conservé à la bibliothèque de l'Institut de France et porte la côte MS-1595.}.} On peut paraphraser cette définition ainsi~:~étant données deux droites $NB, NC$ se coupant en $N$, et $F$ un point pris hors de ces deux droites, si $\tau$ est une droite telle que, pour toute droite passant par $F$ coupant $NB$ en $Y$, $NC$ en $X$ et $\tau$ en $G$ les $F,G;X,Y$ forment une involution, alors $\tau$ est appelée \emph{traversale du point $F$ eu égard à la figure $NB,NC$.} 

En nommant les droites passant par le but $F$ les \emph{ordonnées,} Desargues se place clairement dans le sillage d'Apollonius, affirmant ainsi qu'il en généralise la théorie des diamètres. Cependant, les \textit{Advertissements} montrent qu'il perçoit que la définition qu'il donne ici pose problème et ne se rattache pas vraiment à une conique générale. Il y introduit en particulier la notion d'\textit{ordinale,} qui est une ordonnée qui ne couperait pas la figure, ce qui n'a pas de sens dans le cas où cette figure est $NB,NC$. Il pose également dans ces advertissements les définitions de \textit{point traversal, de diamétrale et de diamétraversale.}


Étant donnée une ordonnée, son point de rencontre avec la traversale correspondante s'appelle son \textit{point traversal.} Si ce point traversal est le milieu du segment d'extrêmités les points de concours d'une ordonnée avec la figure $NB,NC$, cette traversale est nommée une \textit{diamétrale} ou, si l'on veut insister sur son caractère de traversale, une \textit{diamétraversale.} Le point $F$ est alors à distance infinie, les ordonnées sont parallèles et l'on retrouve le cas apollonien des diamètres.

Notons qu'outre généraliser la notion de diamètre, Desargues poursuit aussi un but unificateur en affichant une volonté de traiter de manière unifiée les trois grands types de coniques, comme l'illustre la phrase qu'il demande d'insérer après la ligne 36 de la page 15\footnote{Advertissements.}~:~«~\textit{Par Advis,} Les plus remarquables proprietez des coupes de rouleau, sont communes à toutes especes \& les noms d'\textit{Ellipse, Parabole, \& Hyperbole,} ne leur ont esté donnez qu'à raison d'evenemens qui sont hors d'elles \& de leur nature.~» Nous verrons plus bas ce que l'on peut mettre sous le mot de «~nature~» employé ici par Desargues.

S'ensuivent, aux lignes 45 à 61 de la page 15, des considérations de comparaison de tailles entre les segments donnés sur une ordonnée par l'intersection avec la figure et la traversale, qui peuvent rappeler ce qu'il mentionne dans sa lettre à Mersenne du 4 avril 1638\footnote{\textit{Voir} \cite{taton}, pages 80 à 86.}, mais auxquelles les auteurs du présent article n'ont pas encore réussi à donner sens. 

Il revient, au début de la page 16, sur le fait que la notion de traversale est bien liée à une conique quelconque, même si rien, pour l'instant, ne lui permet d'étayer une telle affirmation. Le reste de la page 16 est dévolu à des considérations sur divers cas de figures, suivant que le but $F$ dont on considère la traversale est à distance finie ou non, qu'il est sur la conique ou non, et suivant que les ordonnées coupent ou ne coupent pas la conique. Les \textit{advertissements} contiennent quelques clarifications intéressantes. Il y énonce en particulier le fait que les points de la traversale du centre d'une conique sont tous à distance infinie, introduisant au passage la \textit{droite à l'infini} du plan. Le tout laisse une impression de grande confusion, confusion qui sera renforcée par ce qui va suivre.

À la fin de la page 16 a lieu l'une de ces transitions brutales dont Desargues est familier~:~il semble mettre de côté la notion de traversale et énonce son \textit{théorème d'involution.} Ce résultat est celui qui lui permettra plus loin de revenir à la traversale en la rattachant à une conique donnée et il va, avec sa démonstration et l'étude de ses cas particuliers, l'occuper jusqu'à la page 20. C'est là, aux lignes 33 et suivantes, qu'il revient à la traversale (d'un point $F$, eu égard à une conique quelconque $\cC$) en justifiant de son existence par une construction. Celle-ci s'effectue bien évidemment\footnote{Rappelons que le théorème d'involution traite des coniques passant par les bornes d'un quadrangle.} en invoquant un quadrangle inscrit dans la conique en question et Desargues laisse entendre, mais n'affirme ni ne prouve vraiment, qu'elle ne dépend pas du quadrangle choisi. Dit autrement, il n'est pas clair à la lecture du \textit{Brouillon} que toute sa construction ne soit pas, en fait, relative au pinceau de toutes les coniques passant par les bornes du quadrangle. Nous verrons cependant que son développement contient un lemme d'incidence qui permet facilement de montrer que la traversale est bien «~eu égard à la coupe de rouleau~» et est indépendante du quadrilatère inscrit pris pour la construire. 

Il détaille ensuite l'involution induite sur une droite $\omega$ passant par $F$, en précisant quels en sont les n{\oe}uds moyens~:~il énonce ici clairement que ce sont le but $F$ des ordonnées et le point d'intersection de $\omega$ avec la traversale de $F$, dit autrement le \emph{point traversal} $H$.  Cette involution est induite par le pinceau des coniques basé sur le quadrangle. Cependant, la conique elle-même induit une autre involution sur la droite $\omega$, dont $F$ et $H$ sont maintenant n{\oe}uds extrêmes, et dont les n{\oe}uds moyens sont les points d'intersection $L$ et $M$ de $\omega$ avec la conique. C'est cette seconde involution qui est induite par la conique elle-même, et non par le pinceau dont elle fait partie.

Jamais Desargues ne précise qu'il obtient ainsi \emph{deux involutions,} et cependant il passe de l'involution induite par le pinceau basé sur le quadrangle à l'involution induite par la conique, ce qui apparaît clairement dans plusieurs passages, dont celui des lignes 30 à 45 de la page 22. Il y montre en outre que la conique induit aussi une involution sur les \emph{ordinales,} c'est-à-dire sur les droites ne coupant pas la conique. Ainsi toute droite peut être considérée comme la traversale d'un certain point, son pôle, et il démontre toutes les propriétés connues de la polarité, comme par exemple qu'elle se comporte correctement vis-à-vis de l'incidence~:~des points alignés ont des traversales concourantes et des droites concourantes sont traversales de points alignés. Enfin, il conclut sur le fait que les propriétés les plus importantes des coniques et les notions qui leurs sont rattachées sont affaire d'incidence et sont donc susceptibles de passer d'un plan de coupe à un autre, donc d'une conique à une autre, par des projections centrales depuis le sommet du cône, éclairant la phrase des \textit{Advertissements} citée plus haut concernant la «~nature~» des coniques.

On peut en conclure que, grâce à théorie de l'involution, Desargues découvre la correspondance points-droites définie par une conique, et qu'à partir de cette correspondance il montre que la conique permet également d'induire, dans \emph{toute} droite du plan\footnote{Non tangente à la conique.} une involution, donnant ainsi une théorie complète de la polarité.

Nous allons maintenant donner une analyse détaillée des arguments donnés par Desargues, en les complétant, et en y remettant de l'ordre, afin de rendre justice à sa magnifique théorie des traversales.

\section{La théorie de la traversale et de la polarité dans le \textit{Brouillon Project}~:~analyse détaillée et mise en forme}
Nous allons dans cette section procéder à une analyse détaillée de la théorie de la traversale et de la polarité développées dans le \textit{Brouillon Project} à partir de la page 20. Mais nous allons également tenter d'en donner une présentation qui la rende complète, en comblant les lacunes laissées par Desargues, et compréhensible au lecteur un tant soit peu averti dans le domaine de la géométrie. Nous respecterons la terminologie et les façons de procéder de Desargues mais nous permettrons de présenter les rapports de grandeurs sous forme de fractions, facilitant ainsi leur manipulation. Enfin nous parlerons de traversale d'un point et non, comme le fait Desargues, des droites de l'ordonnance en ce point. Même si ces deux manières de dire ont des contenus sémantiques quelque peu différents, elles sont, d'un point de vue mathématique, équivalentes.

\subsection{La traversale des droites d'une ordonnance eu égard à une conique}
Nous partirons de la deuxième définition de la traversale d'un point eu égard à une conique, donnée aux lignes 33 et suivantes de la page 20.  Nous pourrons ensuite, en suivant et en complétant les raisonnements de Desargues, démontrer qu'une telle traversale existe, qu'elle est unique au sens où elle est indépendante du procédé utilisé pour la construire, et qu'elle est un objet bien fondé en géométrie tant euclidienne qu'arguésienne, c'est-à-dire projective. 
\begin{definition} Soient $F$ un point et $\cC$ une conique. Une droite $\tau$ est dite \emph{traversale} du point $F$ \emph{eu égard à la conique} $\cC$ si, pour toute droite $\delta$ passant par $F$  coupant $\cC$ (\textit{resp.} $\tau$) en $L$ et en $M$ (\textit{resp.} en $H$), les points $L,M,F,H$ sont en involution.
\end{definition}
Deux problèmes se posent immédiatement à la lecture d'une telle définition. Tout d'abord, une telle traversale existe-t-elle et, si oui, peut-on la construire, ce qui donnerait en sus son unicité ? Ensuite, l'involution dont il est question dans la définition étant de \emph{quatre points seulement,} quels en sont les n{\oe}uds extrêmes et quels en sont les n{\oe}uds moyens doubles ? Cette donnée est indispensable afin de bien définir \emph{une} involution à partir de ces quatre points $L,M,F,H$.

\begin{figure}[!ht]
\centering
\definecolor{qqwuqq}{rgb}{0.,0.39215686274509803,0.}
\definecolor{uuuuuu}{rgb}{0.26666666666666666,0.26666666666666666,0.26666666666666666}
\definecolor{qqqqff}{rgb}{0.,0.,1.}
\begin{tikzpicture}[line cap=round,line join=round,>=triangle 45,x=0.5995203836930456cm,y=0.5995203836930456cm]
\clip(-5.47,-5.63) rectangle (11.21,9.25);
\draw [rotate around={161.3411074719956:(6.718011949781292,-1.614238588797848)},line width=0.8pt] (6.718011949781292,-1.614238588797848) ellipse (2.1686052677999412cm and 1.829789151732657cm);
\draw [line width=0.8pt,domain=-5.47:11.21] plot(\x,{(--0.1124-0.62*\x)/-3.3});
\draw [line width=0.8pt,domain=-5.47:11.21] plot(\x,{(--10.4752--1.36*\x)/-5.48});
\draw [line width=0.8pt,domain=-5.47:11.21] plot(\x,{(--14.1852-3.54*\x)/-0.54});
\draw [line width=0.8pt,domain=-5.47:11.21] plot(\x,{(-43.1888--5.52*\x)/-1.64});
\draw [line width=0.8pt] (4.12,0.74)-- (9.06,-4.16);
\draw [line width=0.8pt] (3.58,-2.8)-- (7.42,1.36);
\draw [line width=0.8pt,color=qqwuqq,domain=-5.47:11.21] plot(\x,{(--50.626487701184566-9.161210824141644*\x)/0.24191301396863363});
\begin{scriptsize}
\draw [fill=qqqqff] (3.58,-2.8) circle (2.5pt);
\draw[color=qqqqff] (3.11,-3.035) node {$E$};
\draw [fill=qqqqff] (4.12,0.74) circle (2.5pt);
\draw[color=qqqqff] (3.77,1.285) node {$D$};
\draw [fill=qqqqff] (7.42,1.36) circle (2.5pt);
\draw[color=qqqqff] (7.64,1.855) node {$C$};
\draw [fill=qqqqff] (9.06,-4.16) circle (2.5pt);
\draw[color=qqqqff] (9.44,-4.595) node {$B$};
\draw [fill=uuuuuu] (-4.305596023130768,-0.8429907679821443) circle (1.5pt);
\draw[color=uuuuuu] (-4.09,-0.395) node {$N$};
\draw [fill=uuuuuu] (5.3020213056541925,8.488806337066372) circle (1.5pt);
\draw[color=uuuuuu] (5.72,8.755) node {$F$};
\draw [fill=uuuuuu] (5.543934319622826,-0.6724044870752725) circle (1.5pt);
\draw[color=uuuuuu] (5.9,-0.605) node {$G$};
\draw [fill=uuuuuu] (5.499792660032231,0.9992337724909043) circle (1.5pt);
\draw[color=uuuuuu] (5.78,0.775) node {$X$};
\draw [fill=uuuuuu] (5.613441988227601,-3.304649836494441) circle (1.5pt);
\draw[color=uuuuuu] (5.9,-2.915) node {$Y$};
\end{scriptsize}
\end{tikzpicture}
\caption{Les points $F,X,G,Y$ sont en involution.}\label{Traversale-01}
\end{figure}
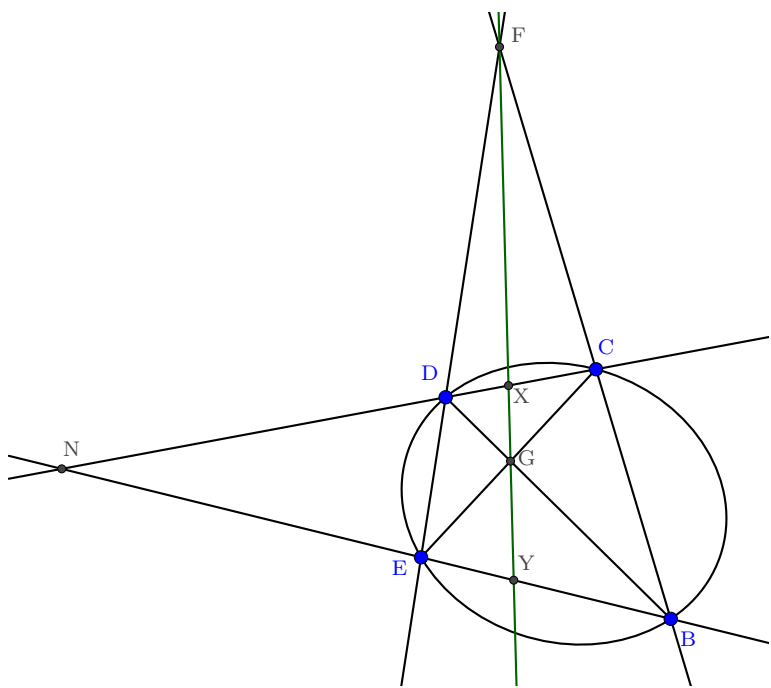

\vs
Supposons pour simplifier les choses que $F$ n'est pas sur la conique. Soient deux droites passant par $F$ et coupant $\cC$ en les points $B,C$ et $D,E$ respectivement, de sorte que $B,C,D,E$ sont les bornes d'un quadrangle inscrit sur $\cC$ dont deux bornales accouplées $BC$ et $DE$ se coupent en $F$. Les deux autres couples de bornales $CD,BE$ et $DB,EC$ se rencontrent respectivement aux points $N$ et $G$. Traçons la droite $FG$~:~elle donne en les droites $NC$ et $NB$ les points $X$ et $Y$, \textit{voir} la figure \ref{Traversale-01}. Rappelons que, pour Desargues, deux droites du même plan sont toujours d'une même ordonnance, et qu'il raisonne avec les points à l'infini comme avec les points ordinaires~:~il considère ses démonstrations valables même si l'un des points $F,N$ ou $G$ est à distance infinie.

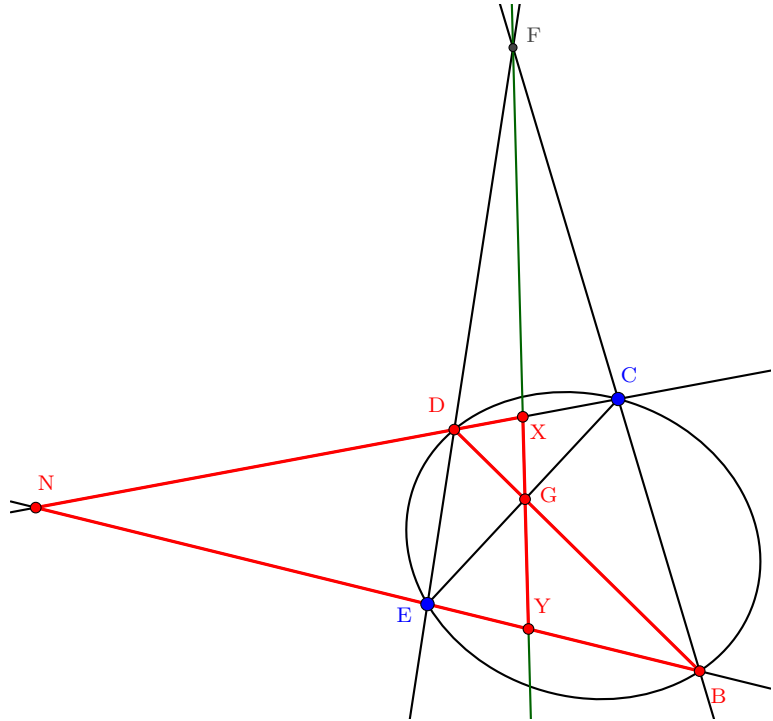
\begin{figure}[!ht]
\centering
\definecolor{qqwuqq}{rgb}{0.,0.39215686274509803,0.}
\definecolor{uuuuuu}{rgb}{0.26666666666666666,0.26666666666666666,0.26666666666666666}
\definecolor{ffqqqq}{rgb}{1.,0.,0.}
\definecolor{qqqqff}{rgb}{0.,0.,1.}
\begin{tikzpicture}[line cap=round,line join=round,>=triangle 45,x=0.6535947712418297cm,y=0.6535947712418297cm]
\clip(-4.81,-5.13) rectangle (10.49,9.36);
\draw [rotate around={161.3411074719956:(6.718011949781292,-1.614238588797848)},line width=0.8pt] (6.718011949781292,-1.614238588797848) ellipse (2.364204958621111cm and 1.9948289575752094cm);
\draw [line width=0.8pt,domain=-4.81:10.49] plot(\x,{(--0.1124-0.62*\x)/-3.3});
\draw [line width=0.8pt,domain=-4.81:10.49] plot(\x,{(--10.4752--1.36*\x)/-5.48});
\draw [line width=0.8pt,domain=-4.81:10.49] plot(\x,{(--14.1852-3.54*\x)/-0.54});
\draw [line width=0.8pt,domain=-4.81:10.49] plot(\x,{(-43.1888--5.52*\x)/-1.64});
\draw [line width=0.8pt] (4.12,0.74)-- (9.06,-4.16);
\draw [line width=0.8pt] (3.58,-2.8)-- (7.42,1.36);
\draw [line width=0.8pt,color=qqwuqq,domain=-4.81:10.49] plot(\x,{(--50.626487701184566-9.161210824141644*\x)/0.24191301396863363});
\draw [line width=1.2pt,color=ffqqqq] (-4.305596023130768,-0.8429907679821443)-- (5.499792660032231,0.9992337724909043);
\draw [line width=1.2pt,color=ffqqqq] (4.12,0.74)-- (9.06,-4.16);
\draw [line width=1.2pt,color=ffqqqq] (5.499792660032231,0.9992337724909043)-- (5.613441988227601,-3.304649836494441);
\draw [line width=1.2pt,color=ffqqqq] (-4.305596023130768,-0.8429907679821443)-- (9.06,-4.16);
\begin{scriptsize}
\draw [fill=qqqqff] (3.58,-2.8) circle (2.5pt);
\draw[color=qqqqff] (3.11,-3.015) node {$E$};
\draw [fill=ffqqqq] (4.12,0.74) circle (2.0pt);
\draw[color=ffqqqq] (3.77,1.245) node {$D$};
\draw [fill=qqqqff] (7.42,1.36) circle (2.5pt);
\draw[color=qqqqff] (7.64,1.845) node {$C$};
\draw [fill=ffqqqq] (9.06,-4.16) circle (2.0pt);
\draw[color=ffqqqq] (9.44,-4.665) node {$B$};
\draw [fill=ffqqqq] (-4.305596023130768,-0.8429907679821443) circle (2.0pt);
\draw[color=ffqqqq] (-4.09,-0.345) node {$N$};
\draw [fill=uuuuuu] (5.3020213056541925,8.488806337066372) circle (1.5pt);
\draw[color=uuuuuu] (5.72,8.745) node {$F$};
\draw [fill=ffqqqq] (5.543934319622826,-0.6724044870752725) circle (2.0pt);
\draw[color=ffqqqq] (6.02,-0.585) node {$G$};
\draw [fill=ffqqqq] (5.499792660032231,0.9992337724909043) circle (2.0pt);
\draw[color=ffqqqq] (5.81,0.705) node {$X$};
\draw [fill=ffqqqq] (5.613441988227601,-3.304649836494441) circle (2.0pt);
\draw[color=ffqqqq] (5.9,-2.835) node {$Y$};
\end{scriptsize}
\end{tikzpicture}
\caption{La figure secteur $N,D,X,G,Y,B$ pour l'identité $GX/GY=(DX/DN)(BN/BY)$.}\label{Traversale-02}
\end{figure}

Desargues commence par démontrer\footnote{p. 20, l. 40 et suivantes.} que les quatre points $F,X,G,Y$ sont en involution, en utilisant le théorème de Ménélaüs dans diverses figures secteurs\footnote{\textit{Voir} l'article \cite{anglade-briend-2} pour un analyse détaillée de l'usage du théorème de Ménélaüs dans le \textit{Brouillon Project.}} illustrées figure \ref{Traversale-02}.


Il obtient de la sorte les identités
\[
\frac{GX}{GY}=\frac{DX}{DN}\frac{BN}{BY}\;\mbox{et}\; \frac{GX}{GY}=\frac{CX}{CN}\frac{EN}{EY}
\]
et
\[
\frac{FX}{FY}=\frac{DX}{DN}\frac{EN}{EY}\;\mbox{et}\; \frac{FX}{FY}=\frac{CX}{CN}\frac{BN}{BY}.
\]
Il en déduit que
\[
\frac{GX}{GY}=\frac{FX}{FY},
\]
identité qui, au vue de la mention «~en concevant que chacune des deux lettres, X \& Y, est doublée~», que l'on trouve à la ligne 39, doit être envisagée comme une égalité entre rapports de carrés~:~
\[
\frac{XG.XG}{YG.YG}=\frac{XF.XF}{YF.YF}.
\]
Si l'on interprète maintenant cette égalité dans le langage arguésien de l'involution, qui s'énonce en termes d'identités entre rapports de rectangles relatifs et rapports obtenus en prenant les rectangles gémeaux correspondants, nous pouvons en déduire que $F,G,X,Y$ forment une involution de quatre points seulement dont $F$ et $G$ sont les n{\oe}uds moyens doubles et $X,Y$ un couple de n{\oe}uds extrêmes, involution dont nous pourrions donc écrire les couplages sous la forme $F,F;G,G;X,Y$. Notons au passage que si l'on omet cette précision sur les n{\oe}uds moyens, on ne fait ici que retrouver des résultats connus depuis Apollonius, comme le fait remarquer Jean de Beaugrand dans ses \textit{Advis Charitables,}  disant que les points $F,G,X,Y$ forment une division harmonique\footnote{On peut citer la proposition 17 du livre III des \textit{Coniques,} \textit{voir} par exemple \cite{apollonius-rashed}, p. 312.}.  

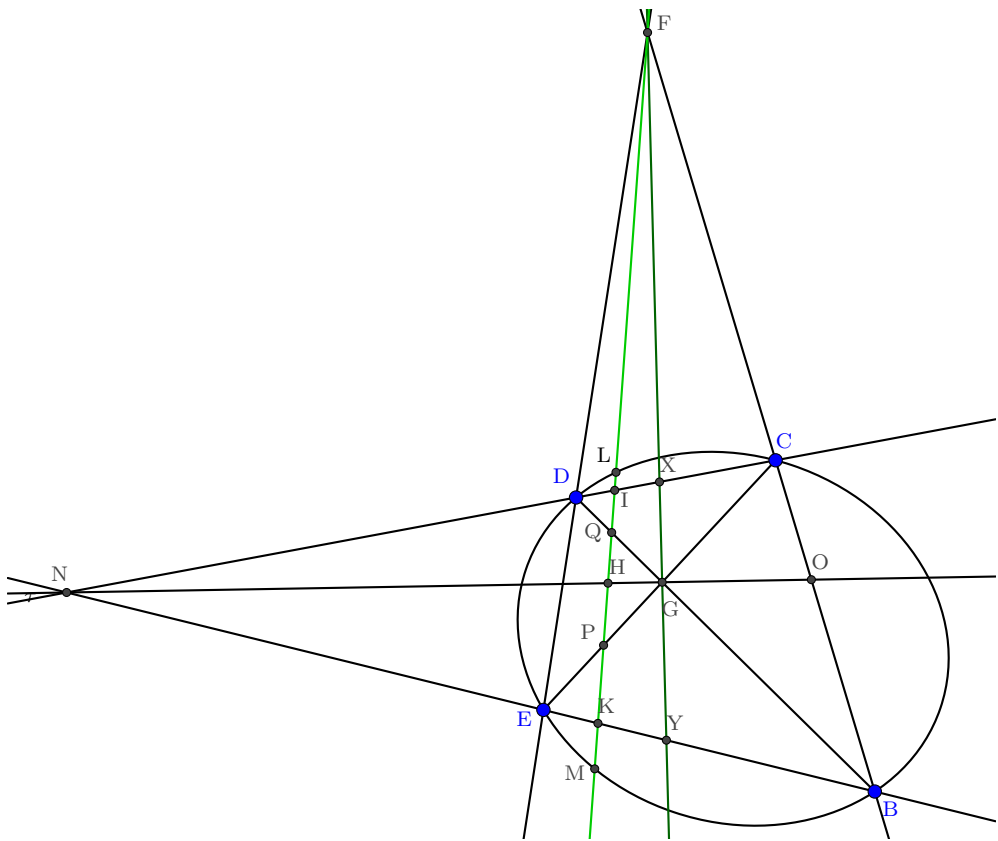
\begin{figure}[!ht]
\centering
\definecolor{qqccqq}{rgb}{0.,0.8,0.}
\definecolor{qqwuqq}{rgb}{0.,0.39215686274509803,0.}
\definecolor{uuuuuu}{rgb}{0.26666666666666666,0.26666666666666666,0.26666666666666666}
\definecolor{qqqqff}{rgb}{0.,0.,1.}
\begin{tikzpicture}[line cap=round,line join=round,x=0.7955936352509176cm,y=0.7955936352509176cm]
\clip(-5.28,-4.94) rectangle (11.06,8.88);
\draw [rotate around={161.3411074719956:(6.718011949781292,-1.614238588797848)},line width=0.8pt] (6.718011949781292,-1.614238588797848) ellipse (2.8778480187866524cm and 2.428222029753422cm);
\draw [line width=0.8pt,domain=-5.28:11.06] plot(\x,{(--0.1124-0.62*\x)/-3.3});
\draw [line width=0.8pt,domain=-5.28:11.06] plot(\x,{(--10.4752--1.36*\x)/-5.48});
\draw [line width=0.8pt,domain=-5.28:11.06] plot(\x,{(--14.1852-3.54*\x)/-0.54});
\draw [line width=0.8pt,domain=-5.28:11.06] plot(\x,{(-43.1888--5.52*\x)/-1.64});
\draw [line width=0.8pt] (4.12,0.74)-- (9.06,-4.16);
\draw [line width=0.8pt] (3.58,-2.8)-- (7.42,1.36);
\draw [line width=0.8pt,color=qqwuqq,domain=-5.28:11.06] plot(\x,{(--50.626487701184566-9.161210824141644*\x)/0.24191301396863363});
\draw [line width=0.8pt,color=qqccqq,domain=-5.28:11.06] plot(\x,{(--34.42614957661839-7.328806337066371*\x)/-0.5220213056541922});
\draw [line width=0.8pt,domain=-5.28:11.06] plot(\x,{(-7.568587535227991--0.1705862809068719*\x)/9.849530342753594});
\begin{scriptsize}
\draw [fill=qqqqff] (3.58,-2.8) circle (2.5pt);
\draw[color=qqqqff] (3.26,-2.95) node {$E$};
\draw [fill=qqqqff] (4.12,0.74) circle (2.5pt);
\draw[color=qqqqff] (3.88,1.11) node {$D$};
\draw [fill=qqqqff] (7.42,1.36) circle (2.5pt);
\draw[color=qqqqff] (7.56,1.69) node {$C$};
\draw [fill=qqqqff] (9.06,-4.16) circle (2.5pt);
\draw[color=qqqqff] (9.32,-4.45) node {$B$};
\draw [fill=black] (4.78,1.16) circle (1.0pt);
\draw[color=black] (4.58,1.45) node {$L$};
\draw [fill=uuuuuu] (-4.305596023130768,-0.8429907679821443) circle (1.5pt);
\draw[color=uuuuuu] (-4.42,-0.53) node {$N$};
\draw [fill=uuuuuu] (5.3020213056541925,8.488806337066372) circle (1.5pt);
\draw[color=uuuuuu] (5.58,8.65) node {$F$};
\draw [fill=uuuuuu] (5.543934319622826,-0.6724044870752725) circle (1.5pt);
\draw[color=uuuuuu] (5.68,-1.11) node {$G$};
\draw [fill=uuuuuu] (5.499792660032231,0.9992337724909043) circle (1.5pt);
\draw[color=uuuuuu] (5.64,1.29) node {$X$};
\draw [fill=uuuuuu] (5.613441988227601,-3.304649836494441) circle (1.5pt);
\draw[color=uuuuuu] (5.76,-3.01) node {$Y$};
\draw [fill=uuuuuu] (4.758630329873502,0.8599850922792641) circle (1.5pt);
\draw[color=uuuuuu] (4.92,0.71) node {$I$};
\draw [fill=uuuuuu] (4.4819897233532515,-3.0238514641898577) circle (1.5pt);
\draw[color=uuuuuu] (4.62,-2.73) node {$K$};
\draw [fill=uuuuuu] (4.78,1.16) circle (1.5pt);
\draw [fill=uuuuuu] (4.427994876718701,-3.781900551263631) circle (1.5pt);
\draw[color=uuuuuu] (4.1,-3.85) node {$M$};
\draw[color=black] (-4.9,-0.92) node {$\tau$};
\draw [fill=uuuuuu] (4.648375436240884,-0.6879148773104405) circle (1.5pt);
\draw[color=uuuuuu] (4.8,-0.41) node {$H$};
\draw [fill=uuuuuu] (8.011135168102461,-0.629674468247313) circle (1.5pt);
\draw[color=uuuuuu] (8.16,-0.33) node {$O$};
\draw [fill=uuuuuu] (4.708504931176668,0.15626029093812388) circle (1.5pt);
\draw[color=uuuuuu] (4.4,0.15) node {$Q$};
\draw [fill=uuuuuu] (4.574688966809071,-1.7224202859568385) circle (1.5pt);
\draw[color=uuuuuu] (4.32,-1.49) node {$P$};
\end{scriptsize}
\end{tikzpicture}
\caption{Diverses involutions~:~$F,I,H,K$, $F,C,O,B$ et $F,Q,P,H$.}\label{Traversale-07}
\end{figure}

Appelons $O$ le point d'intersection de la droite $BC$ avec la droite $GN$ et introduisons une droite quelconque $\delta$ de l'ordonnance en $F$. Elle rencontre $DC$ en $I$, $NG$ en $H$, $EB$ en $K$, $DB$ en $Q$ et $EC$ en $P$, \textit{voir} la figure \ref{Traversale-07} pour une illustation de la situation considérée. 

Nous venons de démontrer que $F,F;G,G;X,Y$ forment une involution. Or cette notion est invariante par projection centrale d'une droite vers une autre, ce qui constitue le \textit{théorème de la ramée} démontré à la page 11 du \textit{Brouillon Project.} D'une ramée en $N$, Desargues peut donc déduire que $F,F;H,H;I,K$, sur $\delta$ et $F,F;O,O;C,B$, sur $BC$, forment des involutions\footnote{p. 20, l. 46-49.}. De cette dernière involution et d'une ramée en $G$ il déduit alors que $F,F;H,H;P,Q$, encore une fois sur $\delta$, forment une involution\footnote{p. 20, l. 50-53.}. 

Il fait alors entrer en scène la conique et note $L$ et $M$ les points d'intersection de celle-ci avec la droite $\delta$. D'après le théorème d'involution sur les pinceaux de coniques, appliqué ici dans le cas particulier où la transversale $\delta$ par laquelle on coupe la figure passe par l'un des points d'intersection (ici $F$) d'un couple de bornales du quadrangle des points bases du pinceau, les trois couples $L,M;I,K;P,Q$ forment une involution\footnote{p. 20, l. 56-61.}. 

Disposant de trois involutions $F,F;H,H;I,K$, ainsi que $F,F;H,H;P,Q$ et $L,M;I,K;P,Q$ sur la droite $\delta$, il utilise l'\textit{agrégativité}\footnote{\textit{Voir} l'article \cite{anglade-briend-1} pour cette notion.} de l'involution pour affirmer finalement que $F,F;H,H;L,M$ forment une involution\footnote{p. 21, l. 1-11.}. 

Nous pouvons en conclure qu'étant donnés la conique $\cC$ et le point $F$, il existe une traversale $\tau$ de $F$ eu égard à $\cC$, à savoir la droite $GN$. Il nous reste à démontrer son unicité et à s'interroger sur la pertinence du choix du couplage $F,F;H,H;L,M$ choisi ci-dessus pour l'involution induite sur $\delta$.
\subsection{Les deux involutions}
L'involution $F,F;H,H;L,M$\footnote{Et ma conique, le H elle aime.} est construite à partir du choix sur la conique $\cC$ de quatre points $B,C,D,E$ qui forment les points bases d'un pinceau de coniques. Cette involution n'est pas, à proprement parler, déduite de la conique particulière $\cC$ qui fait partie de ce pinceau. En effet, si l'on se souvient que par cinq points passe une seule conique et que deux pinceaux différents ne peuvent avoir qu'au plus une conique commune, il est aisé de se convaincre que l'involution donnée dans une droite par le théorème de Desargues sur les pinceaux n'est presque jamais égale à celle induite par \emph{une} conique choisie dans ce pinceau. 



On peut interpréter les lignes 12 et 13 de la page 21 comme suit~:~en y changeant le couplage des points, c'est-à-dire en changeant la manière dont les n{\oe}uds sont considérés comme \textit{correspondant entre eux} ou, si l'on parle en termes d'arbre, en passant à la \textit{souche réciproque,} Desargues change d'involution pour maintenant considérer que $L,L;M,M;F,H$ forment également une involution qui, elle, ne dépendra que de la conique.  Ainsi peut-il écrire que la droite $GN$ est bien traversale des droites de l'ordonnance au but $F$ \emph{eu égard à la conique $\cC$\footnote{p. 21, l. 14-16.}.} Nous nommerons la première involution, celle induite par le pinceau, l'\textit{involution de  Desargues,} car elle est donnée par son théorème d'involution pour les pinceaux, et la seconde l'\textit{involution de polarité,} car nous verrons qu'elle coïncide avec celle qui est donnée par la théorie moderne de la polarité. Notons au passage que dans cette seconde involution, nous perdons les couplages entre les points sur les bornales du quadrangle, comme $I,K$ ou $P$ et $Q$. Nous verrons également plus loin que Desargues a compris que la conique induit une involution de polarité sur la droite $\delta$ même quand celle-ci \emph{ne coupe pas} la conique.
\subsection{La traversale~:~un lemme d'incidence fondamental}
Desargues se lance ensuite un peu précipitamment dans le développement de la théorie de la traversale, avant de s'apercevoir qu'il lui faudrait, pour être complet, démontrer que le résultat de la construction faite ci-dessus ne dépend pas du quadrangle $B,C,D,E$ choisi. Cela est l'objet d'un développement qui court des lignes 27 à 41 de la page 21, même si la preuve qu'il y donne est lacunaire.

Rappelons les données~:~$\cC$ est une conique et $F$ est un point pris hors de celle-ci. D'après ce que nous avons démontré plus haut, ce point admet une traversale $\tau$ eu égard à $\cC$. Soit $N$ un point sur $\tau$ et traçons par $N$ une droite rencontrant $\cC$ en deux points $D$ et $C$. Traçons la droite $FD$ ; outre en $D$, elle recontre la conique en un second point $E$. 

À la ligne 32 de la page 21, Desargues affirme que $NE$ et $FC$ se rencontrent en un point qui est \emph{sur la conique,} ce qui est un \emph{lemme fondamental} pour bien fonder la théorie de la traversale, comme nous allons le voir dans la suite. 

\begin{figure}[!ht]
\centering
\definecolor{qqqqff}{rgb}{0.,0.,1.}
\definecolor{qqqqcc}{rgb}{0.,0.,0.8}
\definecolor{uuuuuu}{rgb}{0.26666666666666666,0.26666666666666666,0.26666666666666666}
\begin{tikzpicture}[line cap=round,line join=round,>=triangle 45,x=0.750341064120055cm,y=0.750341064120055cm]
\clip(-3.86,-7.44) rectangle (10.8,4.14);
\draw [rotate around={20.38684449428258:(5.492200034892476,-3.73506881234851)},line width=1.2pt,color=qqqqcc] (5.492200034892476,-3.73506881234851) ellipse (3.940931330776807cm and 2.419712569274613cm);
\draw [line width=0.8pt,domain=-3.86:10.8] plot(\x,{(--14.2448-6.5*\x)/-2.02});
\draw [line width=1.2pt,color=qqqqcc,domain=-3.86:10.8] plot(\x,{(-61.2056--3.16*\x)/13.8});
\draw [line width=0.8pt,domain=-3.86:10.8] plot(\x,{(-16.1196--2.46*\x)/4.68});
\draw [line width=0.8pt,domain=-3.86:10.8] plot(\x,{(--23.559056188238863-4.044466925682913*\x)/2.651501946261775});
\draw [line width=0.8pt,domain=-3.86:10.8] plot(\x,{(-20.698082878099786-0.19573252790010542*\x)/3.85468004517566});
\begin{scriptsize}
\draw [fill=uuuuuu] (1.34,-2.74) circle (1.5pt);
\draw[color=uuuuuu] (0.96,-2.49) node {$D$};
\draw[color=qqqqcc] (3.,-0.93) node {$c$};
\draw [fill=qqqqff] (3.36,3.76) circle (2.5pt);
\draw[color=qqqqff] (3.02,3.87) node {$F$};
\draw [fill=uuuuuu] (0.5146800451756604,-5.395732527900106) circle (1.5pt);
\draw[color=uuuuuu] (0.16,-5.53) node {$E$};
\draw [fill=qqqqff] (-3.34,-5.2) circle (2.5pt);
\draw[color=qqqqff] (-3.56,-4.75) node {$N$};
\draw[color=qqqqcc] (-1.92,-5.03) node {$\tau$};
\draw [fill=uuuuuu] (6.011501946261775,-0.284466925682914) circle (1.5pt);
\draw[color=uuuuuu] (6.04,0.21) node {$C$};
\draw [fill=uuuuuu] (8.45528915449245,-5.798939921462741) circle (1.5pt);
\draw[color=uuuuuu] (8.3,-5.51) node {$B$};
\draw [fill=uuuuuu] (9.66706057625771,-5.86047112011983) circle (1.5pt);
\draw[color=uuuuuu] (9.9,-5.55) node {$B'$};
\draw [fill=uuuuuu] (9.219468763873362,-5.177737213796593) circle (1.5pt);
\draw[color=uuuuuu] (9.22,-4.57) node {$B''$};
\draw [fill=uuuuuu] (0.8754888329844646,-4.234714151287616) circle (1.5pt);
\draw[color=uuuuuu] (1.1,-4.39) node {$R$};
\draw [fill=uuuuuu] (7.592827514457618,-2.6965409459647773) circle (1.5pt);
\draw[color=uuuuuu] (7.66,-2.31) node {$O$};
\end{scriptsize}
\end{tikzpicture}
\caption{Les points $B,B'$ et $B"$ sont confondus et sur la conique.}\label{Construction-Polaire-00}
\end{figure}
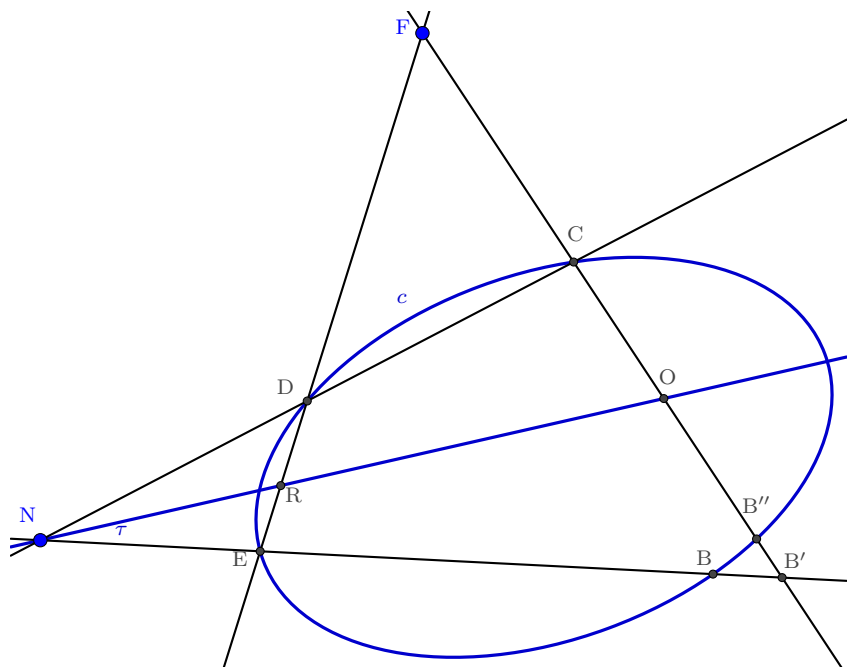

Notons $B$ le point d'intersection de $\cC$ avec $NE$, $B'$ celui de $NE$ avec $FC$ et $B"$ le point d'intersection de $FC$ avec $\cC$. Nommons enfin $O$ et $R$ les points d'intersection respectifs de $FC$ et $FD$ avec $\tau$, \textit{voir} la figure \ref{Construction-Polaire-00} pour une illustration de cette construction. Comme $\tau$ est traversale de $F$ eu égard à $\cC$, les points $C,C;B",B";F,O$ forment une involution, et il en est de même de $D,D;E,E;F,R$. 

\begin{figure}[!ht]
\centering
\definecolor{ffqqtt}{rgb}{1.,0.,0.2}
\definecolor{qqqqff}{rgb}{0.,0.,1.}
\definecolor{qqqqcc}{rgb}{0.,0.,0.8}
\definecolor{uuuuuu}{rgb}{0.26666666666666666,0.26666666666666666,0.26666666666666666}
\begin{tikzpicture}[line cap=round,line join=round,x=0.7777777777777782cm,y=0.7777777777777782cm]
\clip(-5.74,-7.78) rectangle (12.26,4.48);
\draw [rotate around={20.38684449428258:(5.492200034892476,-3.73506881234851)},line width=1.2pt,color=qqqqcc] (5.492200034892476,-3.73506881234851) ellipse (4.0850340723668275cm and 2.508190948070311cm);
\draw [line width=0.8pt,domain=-5.74:12.26] plot(\x,{(--14.2448-6.5*\x)/-2.02});
\draw [line width=1.2pt,color=qqqqcc,domain=-5.74:12.26] plot(\x,{(-61.2056--3.16*\x)/13.8});
\draw [line width=0.8pt,domain=-5.74:12.26] plot(\x,{(-16.1196--2.46*\x)/4.68});
\draw [line width=0.8pt,domain=-5.74:12.26] plot(\x,{(--23.559056188238863-4.044466925682913*\x)/2.651501946261775});
\draw [line width=0.8pt,domain=-5.74:12.26] plot(\x,{(-20.698082878099786-0.19573252790010542*\x)/3.85468004517566});
\draw [line width=0.8pt,domain=-5.74:12.26] plot(\x,{(--47.08090947608595-0.21799531410351225*\x)/-8.704788718697701});
\draw [line width=0.8pt] (1.34,-2.74)-- (9.219468763873362,-5.177737213796593);
\draw [line width=0.8pt] (0.5146800451756604,-5.395732527900106)-- (6.011501946261775,-0.284466925682914);
\draw [line width=0.8pt,color=ffqqtt,domain=-5.74:12.26] plot(\x,{(-28.371487737487698--2.295467848669913*\x)/6.78759567146911});
\begin{scriptsize}
\draw [fill=uuuuuu] (1.34,-2.74) circle (1.5pt);
\draw[color=uuuuuu] (0.96,-2.49) node {$D$};
\draw [fill=qqqqff] (3.36,3.76) circle (2.5pt);
\draw[color=qqqqff] (3.02,3.87) node {$F$};
\draw [fill=uuuuuu] (0.5146800451756604,-5.395732527900106) circle (1.5pt);
\draw[color=uuuuuu] (0.16,-5.53) node {$E$};
\draw [fill=qqqqff] (-3.34,-5.2) circle (2.5pt);
\draw[color=qqqqff] (-3.56,-4.75) node {$N$};
\draw[color=qqqqcc] (-5.26,-5.75) node {$\tau$};
\draw [fill=uuuuuu] (6.011501946261775,-0.284466925682914) circle (1.5pt);
\draw[color=uuuuuu] (6.04,0.21) node {$C$};
\draw [fill=uuuuuu] (8.45528915449245,-5.798939921462741) circle (1.5pt);
\draw[color=uuuuuu] (8.3,-5.51) node {$B$};
\draw [fill=uuuuuu] (9.66706057625771,-5.86047112011983) circle (1.5pt);
\draw[color=uuuuuu] (9.9,-5.55) node {$B'$};
\draw [fill=uuuuuu] (9.219468763873362,-5.177737213796593) circle (1.5pt);
\draw[color=uuuuuu] (9.22,-4.57) node {$B''$};
\draw [fill=uuuuuu] (0.8754888329844646,-4.234714151287616) circle (1.5pt);
\draw[color=uuuuuu] (1.1,-4.39) node {$R$};
\draw [fill=uuuuuu] (7.592827514457618,-2.6965409459647773) circle (1.5pt);
\draw[color=uuuuuu] (7.66,-2.31) node {$O$};
\draw [fill=uuuuuu] (-3.923833581389661,-5.506886882525334) circle (1.5pt);
\draw[color=uuuuuu] (-3.84,-5.73) node {$N'$};
\draw [fill=uuuuuu] (2.863762090079449,-3.2114190338554214) circle (1.5pt);
\draw[color=uuuuuu] (2.82,-2.81) node {$G'$};
\draw[color=ffqqtt] (-5.28,-6.07) node {$\tau'$};
\draw [fill=uuuuuu] (0.99734086879691,-3.8426160162475687) circle (1.5pt);
\draw[color=uuuuuu] (1.34,-3.47) node {$R'$};
\draw [fill=uuuuuu] (7.010909716554203,-1.8089129770678485) circle (1.5pt);
\draw[color=uuuuuu] (7.22,-1.51) node {$O'$};
\end{scriptsize}
\end{tikzpicture}
\caption{Les deux traversales $\tau$ et $\tau'$ sont confondues.}\label{Construction-Polaire-04}
\end{figure}
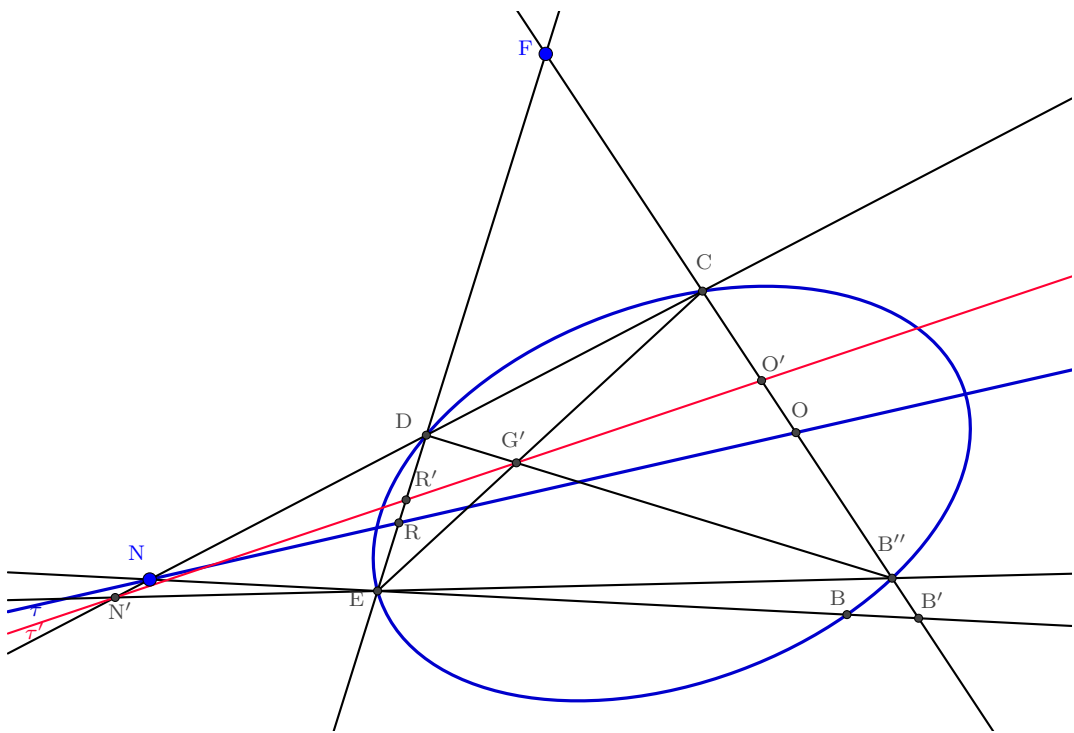

Par ailleurs les points $B",C,D,E$ forment un quadrangle inscrit sur la conique $\cC$. Notons $F,N',G'$ les points d'intersection de ses couples de bornales, \textit{voir} la figure \ref{Construction-Polaire-04}.

D'après ce qui précède, la droite $\tau'=G'N'$ est elle aussi traversale de $F$ eu égard à la conique $\cC$. Notons $O'$ et $R'$ les points d'intersection de $\tau'$ avec $FB"$ et $FD$. Alors $C,C;B",B"; F,O'$ et $D,D; E,E; F,R'$ forment des involutions de sorte que, par une propriété générale des involutions (en fait, des divisions harmoniques) on a nécessairement $O'=O$ et $R'=R$, de sorte que $\tau'=\tau$. Ainsi $G'$ et $N'$ sont sur $\tau$ et comme $N$ est point d'intersection de $DC$ avec $\tau$, et que $N'$ est point d'intersection de $DC$ avec $\tau'$ qui n'est autre que $\tau$ on en déduit que $N'=N$. Il en découle, du fait qu'en outre les trois droites $B"E, N'E$ et $NE$ sont confondues, que $B=B"$ puis enfin que $B=B'=B"$, ce qu'il fallait démontrer.  

L'argument donné par Desargues dans le \textit{Brouillon Project} est incomplet. En effet, il y affirme de but en blanc que $\tau$ est aussi traversale de $F$ eu égard à la figure $NB,NC$, ce qui, entraînant que $C,C;B,B;F,O$ forment aussi une involution implique que $B=B'$. Il pèche ici par excès d'optimisme~:~l'affirmation que $\tau$ est traversale de $F$ eu égard à la figure $NB,NC$ est vraie mais doit être démontrée. Il lui faudrait pour cela prouver, par exemple, que $B'E$ et $CD$ se coupent en un point situé sur $\tau$, comme nous venons de le démontrer. L'affirmation énoncée ci-dessus permet de donner un procédé de génération de la conique par incidence, ce qui n'est pas sans rappeler ce que fait Philippe de la Hire dans ses \textit{Planiconiques\footnote{Les \textit{Planiconiques} sont un complément à l'ouvrage \cite{planiconiques}. On pourra se reporter à la thèse de \v{S}\'{\i}r \cite{sir}.}}. 

Notons en outre que ce qui précède permet de justifier une affirmation qui se trouve un peu avant ce développement dans le texte de Desargues, p. 21 l. 16 et 17, où, en reprenant les notations qui précèdent, il est dit que $FN$ est traversale de $G$ et $GF$ traversale de $N$, ce qui permettra plus loin de mettre au jour le caractère symétrique de la correspondance associant à un point sa traversale.

Enfin, il faut remarquer que c'est en considérant la bonne involution, à savoir celle de polarité, que nous avons été en mesure de compléter la preuve de Desargues. 
\subsection{La traversale est «~donnée de position~»}
Desargues affirme, aux lignes 41 à 46, que si l'on se donne une conique $\cC$ et un point $F$ hors de $\cC$, un quadrangle \emph{quelconque} $B,C,D,E$ inscrit sur $\cC$ tel que $BC$ et $DE$ se coupent en $F$, alors la droite joignant les points d'intersection $G$ et $N$ des deux autres couples de bornales est \emph{la} traversale de $F$ eu égard à $\cC$. Desargues ne démontre pas cette affirmation mais nous allons voir qu'on peut la déduire très facilement de ce qui précède. Considérons en effet deux autres points $L,M$ sur $\cC$ de sorte que $LM$ passe par $F$, donnant ainsi naissance à deux autres quadrangles $L,M,B,C$ et $L,M,D,E$ et donc, potentiellement, à deux autres traversales.

Notons $\tau$ la droite $GN$ et $H$ le point de rencontre entre $FL$ et $\tau$. D'après ce que l'on a démontré plus haut, $\tau$ est traversale de $F$ eu égard à $\cC$ et donc les points $L,L;M,M;F,H$ forment une involution. Notons $R$ et $S$ les points de concours respectifs des deux couples de bornales $LB,CM$ et $LC, BM$ du quadrangle $L,M,C,B$. Toujours d'après ce que l'on a démontré plus haut, $RS$ est traversale aussi de $F$ eu égard à $\cC$ et si l'on note $H'$ le point de rencontre de $RS$ et $FL$ on obtient alors que $L,L;M,M;F,H'$ forment une involution, ce qui entraine que $H=H'$ d'où l'on tire aisément, quitte à répéter le même argument avec une autre droite passant par $F$, que $RS=\tau$, achevant de montrer que la construction de la traversale grâce à un quadrangle inscrit est valide et produit toujours la même droite, \emph{la} traversale de $F$ eu égard à $\cC$.

Nous pouvons donc conclure, comme le fait Desargues\footnote{p. 21, l. 53}, qu'«~estant au plan d'une coupe de rouleau donné de position le but $F$, d'une quelconque ordonnance de droictes $FH, FG$, leur commune traversale $GN$, y est aussi donnée de position.~»

\begin{figure}[!ht]
\centering
\includegraphics[width=10cm]{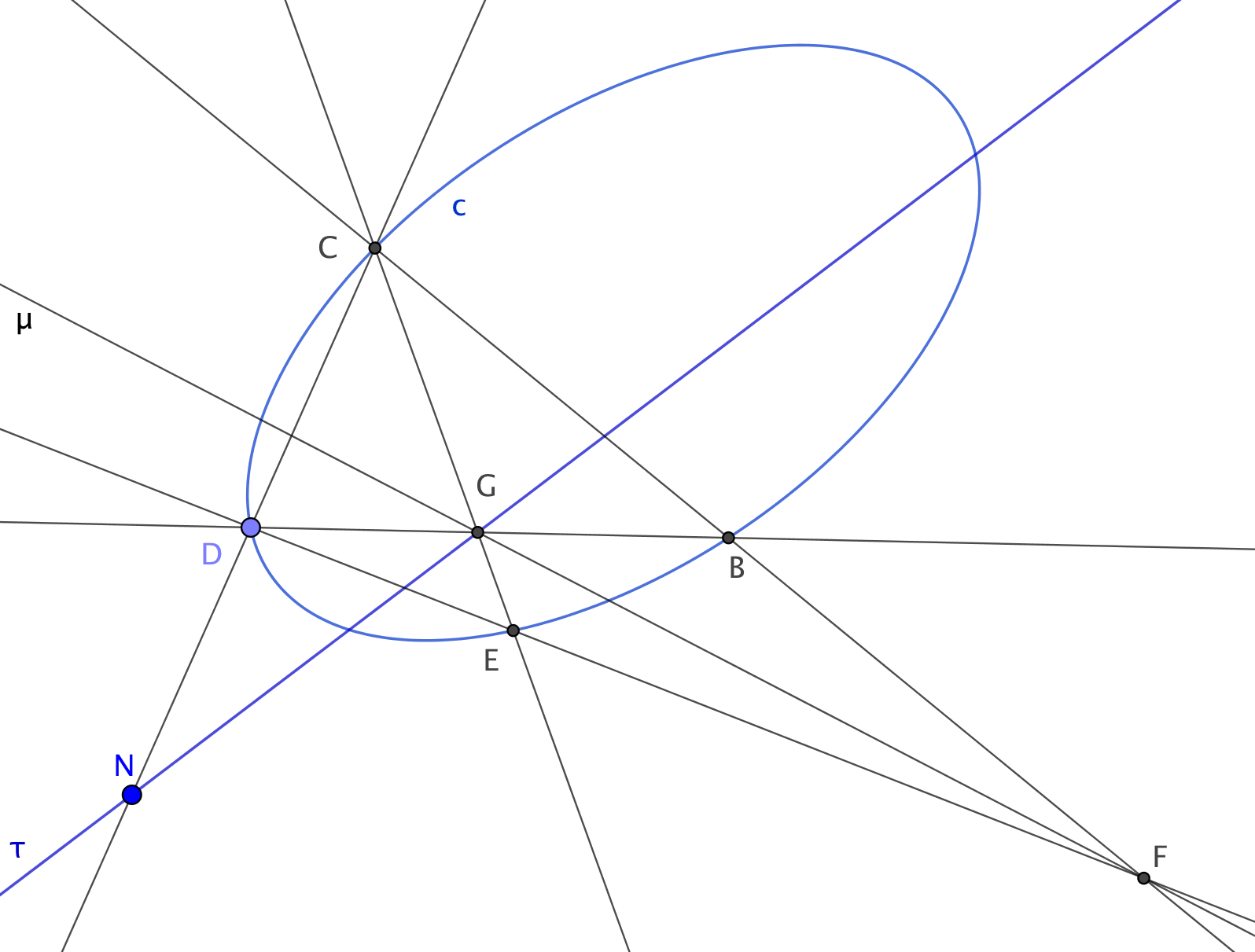}
\caption{Construction du pôle de la droite $\tau$.}\label{Construction-Pole-03}
\end{figure}

La phrase qui suit est en quelque sorte un énoncé dual de celui qui le précède~:~«~Et qu'y estant donnée de position une quelconque traversale ou droicte $GN$, le but de ses ordonnées $F$, y est aussi donné de position.~» Ainsi toute droite peut être considérée comme une traversale, eu égard à la conique $\cC$, d'un point (ou but) que l'on peut construire. C'est une affirmation que l'on rencontre un peu avant dans le texte, à la ligne 23 de la page 21. Voyons comment on peut la déduire facilement du lemme fondamental démontré dans la section précédente, \textit{voir} la figure \ref{Construction-Pole-03}. Donnons-nous la conique $\cC$ et une droite $\tau$. Prenons un point $N$ sur $\tau$ et construisons la traversale $\mu$ de $N$ eu égard à $\cC$, ce qui est possible d'après tout ce qui précède. Cette traversale coupe $\tau$ en un point que nous nommerons $G$, qui est le point traversal de $N$ eu égard à $\cC$ sur la droite $\tau$. Traçons par $N$ une droite qui coupe $\cC$ en deux points $C$ et $D$ et traçons les deux droites $CG,DG$, qui, outre en $C$ et $D$, vont couper $\cC$ en les points respectifs $E$ et $B$. Par un raisonnement très analogue à celui effectué pour démontrer le lemme d'incidence fondamental évoqué plus haut, on démontre que les points $N,E$ et $B$ sont alignés, nous donnant ainsi un quadrangle de construction de la traversale $\mu$ du point $N$. Par unicité d'icelle, le point d'intersection $F$ de $CB$ et $DE$ est sur la droite $\mu$ et l'on en conlut que $\tau=GN$ est traversale de $F$, ce qu'il fallait démontrer. En termes modernes, le point $F$ n'est autre que le \emph{pôle} de la droite $\tau$, qui est sa \emph{polaire.}

\subsection{Polarité~:~la dualité points-droites}
Nous pouvons résumer le résultat des lignes qui précèdent en deux théorèmes, dont le premier exprime l'existence, l'unicité et un moyen de construction de la traversale d'un point eu égard à une conique.
\begin{theoreme}\label{theo-traversale} Soient $\cC$ une conique et $F$ un point hors de $\cC$. Alors $F$ admet, eu égard à la conique $\cC$, une et une seule traversale $T_{F/\cC}$ qui est une droite ayant la propriété suivante~:~pour tout droite $\omega$ passant par $F$ et coupant $\cC$, si l'on note $L,M$ les points d'intersections de $\omega$ avec $\cC$ et $O_\omega$ le point d'intersection de $\omega$ avec $T_{F/\cC}$, alors les points $L,L;M,M;F,O_\omega$ sont en involution. De plus, pour tout quadrilatère $B,C,D,E$ inscrit sur la conique et dont les couples de bornales $ED,BC;CD,BE;DB,EC$ se coupent respectivement en $F$, $N$ et $G$, la droite $GN$ est égale à la traversale $T_{F/\cC}$ de $F$ eu égard à $\cC$. 
\end{theoreme}
Les droites telles $\omega$ passant par $F$ sont appelées les \emph{ordonnées} de la traversale $T_{F/\cC}$. Le point $O_\omega$ est le \emph{point traversal} de l'ordonnée $\omega$. 

Le second théorème exprime la possibilité d'une construction duale, partant d'une droite pour construire un point dont elle est traversale.
\begin{theoreme}\label{theo-pole} Soient $\cC$ une conique et $\tau$ une droite du même plan que celui de $\cC$, non tangente à $\cC$. Alors $\tau$ est traversale d'un point $P_{\tau/\cC}$ et un seul. Si $G$ et $N$ sont deux points pris sur $\tau$ et si l'on construit (grâce à un quadrangle, par exemple) les traversales $T_{G/\cC}$ et $T_{N/\cC}$, alors $P_{\tau/\cC}$ est le point de rencontre de ces deux traversales.
\end{theoreme}
Le point $P_{\tau/\cC}$ s'appelle le \textit{pôle} de la droite $\tau$ eu égard à la conique $\cC$. Desargues ne nomme pas ce point dans le \textit{Brouillon Project} mais, par commodité, nous adopterons pour icelui la terminologie moderne.

Ainsi, toute conique $\cC$ donne naissance à une correspondance points-droites, dont il est de plus aisé, au vu de ce qui précède, de voir qu'elle se comporte de manière satisfaisante eu égard aux propriétés d'incidence. Ainsi Desargues écrit-il\footnote{p. 21, l. 47} que si trois droites sont concourantes, alors leurs pôles sont alignés. De manière duale, ou «~par converse~», comme il l'écrit\footnote{p. 21, l. 50}, si trois points sont alignés, alors leurs traversales sont concourantes. Nous obtenons ce que l'on appelle aujourd'hui une \emph{dualité} projective entre un plan projectif et son dual.

Desargues examine ensuite quelques cas particuliers qu'il a laissés de côté. Il le fait un peu rapidement mais ce qu'il écrit est tout à fait convaincant au vu de ce qui vient d'être démontré.

\subsection{La traversale~:~tangentes, diamétrales, centre}
Au moment de la rédaction du \textit{Brouillon Project,} on ne dispose pas encore d'une théorie vraiment satisfaisante des tangentes à une courbe et l'on doit parfois se contenter de raisonnements \emph{par continuité} qui manquent de la rigueur qu'il acquerront plus tard avec l'introduction de la géométrie analytique. Cette approche, que l'on pourrait qualifier de sensible, est celle que nous adopterons ici.

Commençons par un retour en arrière dans le \textit{Brouillon,} à la page 16 et aux tous premiers développements de la notion de traversale. Desargues y note qu'il existe plusieurs configurations possibles, selon que la traversale coupe ou ne coupe pas la conique, ou selon que le but de ses ordonnées est à l'extérieur ou à l'intérieur de celle-ci. 

Si une traversale $\tau$ rencontre la conique en deux points, alors le but de ses ordonnées est situé à l'extérieur du domaine enserré par cette conique et les ordonnées peuvent donc soit la couper, soit lui être tangente, soit ne pas la rencontrer du tout\footnote{p. 16, l. 36.}. Dans ces deux derniers cas, les ordonnées correspondantes sont nommées des \textit{ordinales\footnote{\textit{Advertissements,} lignes à rajouter après p. 15, l. 44.}.} Réciproquement, si un point $F$ est situé à l'extérieur de la conique, alors sa traversale coupe la conique en deux points. 

Si en revanche une traversale $\tau$ ne rencontre pas la conique, alors le but de ses ordonnées est à l'intérieur de la conique et toutes les ordonnées coupent cette dernière\footnote{p. 16, l. 44.}. Réciproquement, si un point est à l'intérieur de la conique, sa traversale ne la rencontre pas. 


Qu'en est-il maintenant des points situés sur la conique, ou des droites qui ne rencontrent la conique qu'en un seul point ?

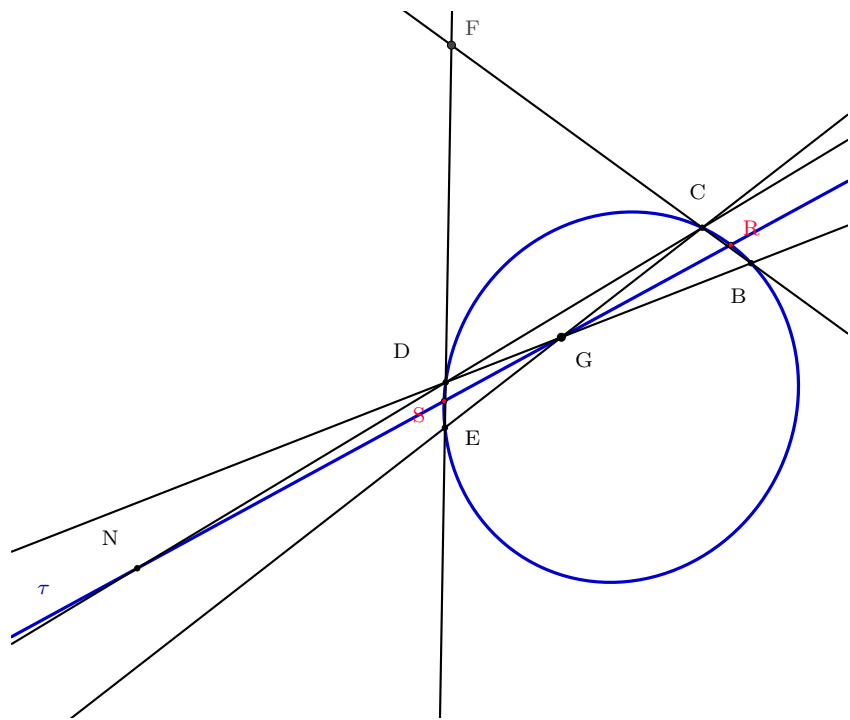
\begin{figure}[!ht]
\centering
\definecolor{ffqqtt}{rgb}{1.,0.,0.2}
\definecolor{uuuuuu}{rgb}{0.26666666666666666,0.26666666666666666,0.26666666666666666}
\definecolor{qqqqcc}{rgb}{0.,0.,0.8}
\begin{tikzpicture}[line cap=round,line join=round,>=triangle 45,x=0.6267806267806254cm,y=0.6267806267806254cm]
\clip(-4.695,-8.285) rectangle (12.855,6.655);
\draw [line width=1.2pt,color=qqqqcc,domain=-4.695:12.855] plot(\x,{(-49.7048--6.84*\x)/12.46});
\draw [rotate around={65.11963314081817:(8.1,-1.5)},line width=1.2pt,color=qqqqcc] (8.1,-1.5) ellipse (2.4875079769155586cm and 2.2975273416030864cm);
\draw [line width=0.8pt,domain=-4.695:12.855] plot(\x,{(-45.89342952891502--7.2049824652356875*\x)/11.862440118613385});
\draw [line width=0.8pt,domain=-4.695:12.855] plot(\x,{(-18.73912781723099--2.520550282134062*\x)/6.411465203810348});
\draw [line width=0.8pt,domain=-4.695:12.855] plot(\x,{(--16.54613684411572-2.3173171128163674*\x)/-2.9588859824811085});
\draw [line width=0.8pt,domain=-4.695:12.855] plot(\x,{(-4.256761561507514--0.9596479354321894*\x)/0.016205203207674046});
\draw [line width=0.8pt,domain=-4.695:12.855] plot(\x,{(-9.500431145891305--0.7512284053666998*\x)/-1.0247374846105917});
\begin{scriptsize}
\draw [fill=black] (-2.06,-5.12) circle (1.0pt);
\draw[color=black] (-2.625,-4.4825) node {$N$};
\draw[color=qqqqcc] (-4.02,-5.5625) node {$\tau$};
\draw [fill=black] (9.802440118613385,2.0849824652356874) circle (1.0pt);
\draw[color=black] (9.705,2.8525) node {$C$};
\draw [fill=black] (4.4157123994136285,-1.186796222265074) circle (1.0pt);
\draw[color=black] (3.495,-0.5225) node {$D$};
\draw [fill=black] (10.827177603223976,1.3337540598689877) circle (1.0pt);
\draw[color=black] (10.56,0.6475) node {$B$};
\draw [fill=black] (6.843554136132276,-0.2323346475806802) circle (1.5pt);
\draw[color=black] (7.32,-0.7025) node {$G$};
\draw [fill=black] (4.399507196205954,-2.1464441576972635) circle (1.0pt);
\draw[color=black] (4.98,-2.3675) node {$E$};
\draw [fill=uuuuuu] (4.5361554537876945,5.945661605383946) circle (1.5pt);
\draw[color=uuuuuu] (4.98,6.3175) node {$F$};
\draw [fill=ffqqtt] (10.4,1.72) circle (1.0pt);
\draw[color=ffqqtt] (10.83,2.0875) node {$R$};
\draw [fill=ffqqtt] (4.381815416359294,-1.583722516220099) circle (1.0pt);
\draw[color=ffqqtt] (3.855,-1.8725) node {$S$};
\end{scriptsize}
\end{tikzpicture}
\caption{La traversale de $F$ relie les points de contact des tangentes à $\cC$ issues de $F$.}\label{ConstructionPolaire-05}
\end{figure}

Une droite qui rencontre une conique en un seul point lui est tangente.  Pour voir le lien entre tangente et  traversale, donnons nous une conique $\cC$ et un point $F$ à l'extérieur du domaine qu'elle enserre. La traversale $\tau$ de $F$ eu égard à $\cC$ est alors une droite qui coupe la conique en deux points $R$ et $S$. Si l'on considère un quadrilatère $B,C,D,E$ inscrit dans $\cC$ et dont deux bornales couplées $CB$ et $DE$ sont de même ordonnance en $F$, nous avons vu que les deux autres couples de bornales se coupent en des points qui sont sur la traversale $\tau$. Lorsque les deux droites $BC$ et $DE$ s'écartent et tendent vers les deux touchantes à $\cC$ passant par $F$, les points $B,C$ tendent à se confondre en $R$ et les points $D,E$ à se confondre en $S$. Ainsi peut-on en conclure que la traversale de $F$ eu égard à $\cC$, lorsque $F$ est extérieur à la conique, est la droite reliant les points de contact des deux tangentes à $\cC$ menées depuis $F$. Cet énoncé est donné par Desargues dans les \textit{Advertissements,} comme devant venir s'insérer après la ligne 5 de la page 22 et la figure \ref{ConstructionPolaire-05} illustre son propos.

En anticipant quelque peu sur la section qui suit, nous pouvons en déduire l'énoncé suivant, que Desargues donne dans le corps du texte\footnote{p. 22, l. 2 à 5.} dans le cas d'une diamétrale, et qu'il généralise dans les \textit{Advertissements} à toute traversale. Prenons sur la conique $\cC$ un point $A$ et menons de ce point $A$ l'une des ordonnées $\omega$ de la traversale $\tau$. Ces deux droites se coupent en un point $H$. Si $R,S$ sont les points d'intersection de $\cC$ et de $\tau$, on peut trouver un point $I$ tel que $R,R;S,S;H,I$ forment une involution (le point $I$ n'est autre que le conjugué harmonique de $H$ par rapport au couple $R,S$). La droite $AI$ est alors tangente à la conique.

Il s'ensuit enfin\footnote{p. 16, l. 33.} que si le point $F$ est \emph{sur} la conique, à distance finie ou infinie, alors sa traversale est «~du corps mesme des ordonnées~», c'est-à-dire passe par $F$, et de plus est tangente en ce point à la conique.

Comme nous l'avons vu à de nombreuses reprises, Desargues adopte sciemment dans le \textit{Brouillon Project} un point de vue projectif, en ce qu'il ne fait pas de distinction, dans ses raisonnements ou ses énoncés, entre les points à distance finie et ceux qui sont à distance infinie. Cependant, du fait qu'une involution de quatre points (c'est-à-dire une division harmonique) donne, quand l'un d'eux est à distance infinie, un segment coupé en deux par son milieu, il lui arrive souvent de préciser ce qu'il advient de ses propositions lorsque l'un des points se trouve à l'infini. Cela lui permet au surplus de montrer que sa théorie généralise celle des diamètres d'Apollonius.

Si $F$ est un point à distance infinie, sa traversale eu égard à une conique $\cC$ est un diamètre au sens d'Apollonius. Afin de montrer que les diamètres ne sont qu'un cas particulier de traversale, Desargues les nomme des \emph{diamétrales\footnote{\textit{Advertissements} de la ligne 44 de la page 15.}} et, s'il les considère comme des traversales d'une ordonnance, il les nomme des \textit{diamétraversales.} Ainsi une diamétraversale est une traversale dont les ordonnées «~mypartissent~» la figure. Dit autrement, soit $\tau$ une diamétraversale et $\omega$ l'une de ses ordonnées rencontrant $\tau$ en $H$ et $\cC$ en $B$ \& $C$ ; alors $H$ est le milieu de $BC$. Réciproquement, si une quelconque ordonnée d'une traversale $\tau$ mypartit la figure, alors le but de ses ordonnées est à distance infinie et $\tau$ est une diamétraversale.

Desargues va pousser le bouchon projectif encore plus loin en introduisant la \textit{droite à l'infini} du plan, qu'il définit comme la traversale du centre d'une conique. Tout d'abord gêné, il écrit,   à la ligne 34 de la page 22~:~«~Quand une traversale est à distance infinie, tout en est inimaginable.~» On peut cependant penser qu'il ait alors fait le raisonnement qui suit~:~si deux diamétraversales $\tau$ et $\tau'$ se rencontrent en un point $O$, alors $O$ doit avoir pour traversale une droite, celle qui relie les deux but des ordonnées $F,F'$ de $\tau$ et $\tau'$. Ces deux buts sont à distance infinie et donc la traversale de $O$ est une droite dont tous les points sont à distance infinie. Le point $O$ est le \emph{centre} de la conique et sa traversale est la droite à l'infinie du plan affine contenant la conique. Réciproquement, le but des ordonnées, ou pôle, de la droite à l'infini doit être sur tous les diamètres de la conique et donc en être le centre. 

Un souci de cohérence avec son affirmation que tout point admet une traversale va donc inciter Desargues à introduire une droite dont tous les points seraient à distance infinie. Après avoir écrit la phrase citée ci-dessus, il rajoute dans les \textit{Advertissements} la note suivante, qu'il voudrait voir insérée après la ligne 3 de la page 16~:~«~Par forme d'éclaircissemens, Quand en un plan, aucun des poincts d'une droicte n'y est à distance finie, cette droicte y est à distance infinie. D'autant qu'en un plan le poinct nommé centre d'une coupe de rouleau n'est qu'un cas d'entre les innombrables buts d'ordonnances de droictes, il ne doit estre icy iamais parlé de centre de coupe de rouleau.~» 

Il précise immédiatement son propos concernant le centre~:~«~D'autant que toute droicte qui passe au sommet d'un rouleau \& au quelconque but d'une ordonnance de droictes au plan de sa base, a une proprieté commune avec celle qui passe au but des diamétrales de la base de ce rouleau, iamais il ne doit estre icy parlé d'essieu de rouleau.~» Nous pouvons interpréter cela comme l'affirmation renouvelée que les propriétés intéressantes des coniques sont invariantes par projection depuis le sommet du cône, c'est-à-dire sont des notions projectives. Au passage, Desargues règle ici définitivement son compte à la méthode du «~triangle par l'axe~» héritée d'Apollonius. C'est la condition \textit{sine qua non} à l'unification de la théorie des coniques et nous en trouverons une illustration dans le traitement qu'il donne des asymptotes à une hyperbole.

\subsection{Polarité~:~l'involution polaire, les ordinales}
Au début de la construction de la traversale, on introduit une involution sur chaque ordonnée qui n'est pas relative à la conique considérée, mais relative au (pinceau des coniques passant par les quatre bornes du) quadrilatère utilisé pour la construire. Rappelons brièvement ce dont il s'agit~:~on se donne une conique $\cC$ et un point $F$ hors de celle-ci. Pour n'importe quel quadrilatère $B,C,D,E$ inscrit sur $\cC$ de sorte que $BC$ et $DE$ se coupent en $F$, les intersections $G$ et $N$ des deux autres couples de bornales permettent de former la droite $GN$ qui est traversale de $F$ eu égard à la conique $\cC$, indépendamment du quadrilatère choisi. Dans ce cas particulier, le théorème d'involution fournit sur toute droite passant par $F$ une involution hyperbolique, dite de Desargues, dont les n{\oe}uds moyens doubles sont situés sur les bornales et, si la traversale ne dépend pas du quadrilatère, l'involution ainsi construite en dépend. 

En passant aux souches réciproques comme nous l'avons décrit plus haut, nous obtenons une autre involution, dite de polarité, qui cette fois-ci est induite par la conique et ne dépend plus du quadrilatère choisi. Nous perdons alors l'accouplement entre les points sur les bornales, qui est contingent à la construction. Comme nous l'avons montré un peu plus haut dans cet article, considérer l'involution de polarité et non celle de Desargues est indispensable pour mener à bien la théorie de la traversale.  

Cependant, un aspect reste encore insatisfaisant concernant la polarité. Dans le théorème d'involution, on induit une involution sur \emph{toute droite} passant par $F$, alors que pour l'instant nous n'avons été en mesure de construire l'involution de polarité que sur les droites coupant la conique. Desargues comble cette lacune à la page 22 du \textit{Brouillon Project,} achevant ainsi l'élaboration de la théorie de la polarité.

Aux lignes 6 à 12 de cette même page, Desargues montre comment construire l'involution de polarité sur une droite $\tau$ qui, hormis n'être pas tangente à la conique $\cC$, est quelconque. Soient $N,Z$ et $A$ trois points sur $\tau$. Les traversales respectives $\mu,\zeta$ et $\alpha$ de ces trois points coupent la droite $\tau$ aux points $G,H$ et $R$. Alors $N,G;Z,H;A,R$ sont trois couples de points d'une involution. Les trois points $N,Z$ et $A$ étant alignés, leurs trois traversales sont concourantes en un point dont $\tau$ est la traversale. Si $\tau$ coupe la conique en deux points $S$ et $T$, on reconnaît dans l'involution $N,G;Z,H;A,R$ l'involution de polarité dont $S$ et $T$ sont n{\oe}uds moyens doubles. Elle est donc une «~involution d'un arbre dont la souche est conséquemment donnée de position~». Si $\tau$ ne coupe pas la conique, l'involution ainsi construite est sans n{\oe}ud moyen double. En effet, les couples formant l'involution, de démêlés qu'ils sont entre eux quand $\tau$ coupe la conique, deviennent mêlés lorsque celle-ci ne la coupe plus. 

Il ressort de ces considérations que les droites d'une ordonnance au but $F$ viennent en deux espèces, eu égard à l'involution de polarité induite sur celles-ci par la conique. Si elles coupent la conique, l'involution est à n{\oe}uds moyens doubles (ou l'arbre est d'espèce à souche dégagée) et on les appelle des \emph{ordonnées.} Si elles ne la coupent pas, l'involution est sans n{\oe}ud moyen double (ou l'arbre est d'espèce à souche engagée) et Desargues les nomme des \emph{ordinales.} Ce mot n'est défini que dans les \textit{Advertissements,} pouvant laisser supposer que l'auteur n'avait pas vu cette subtilité lors de la rédaction du premier jet. 

Une remarque s'impose ici~:~dans la théorie antique, une conique induit des divisions harmoniques sur des droites qui doivent nécessairement la couper. Sur une droite qui ne rencontre pas la conique, la notion de division harmonique n'est pas pertinente mais, en plus, on voit mal comment y définir une quelconque configuration de points. Grâce à sa notion d'involution, Desargues montre que la conique munit toute droite du plan, même celles qui sont \textit{à distance de la conique,} d'une structure géométrique caractéristique. Ce point seul mérite l'introduction de la notion d'involution et c'est ce que les contemporains de Desargues, comme Jean de Beaugrand, ou ses suivants immédiats comme Philippe de la Hire, n'ont pas compris.

Notons également qu'aux lignes suivant immédiatement le passage évoqué ci-dessus\footnote{p. 22, l. 13 à 21.}, Desargues mentionne qu'un certain «~Monsieur Pujoz~» aurait une démonstration générale et dans le plan de cette proposition, et que par ailleurs, étant évidente pour le cercle du fait de la perpendicularité des diamètres, elle se généralise à toute coupe en «~rétablissant le rouleau sur cette coupe~», c'est à dire en utilisant une ramée depuis le sommet du cône. 

Que se passe-t-il lorsque la droite considérée devient tangente à la conique ? Desargues traite de ce cas aux lignes 30 et suivantes de la page 22. Pour simplifier notre propos, revenons aux droites ordonnées au but $F$ et faisons les varier.

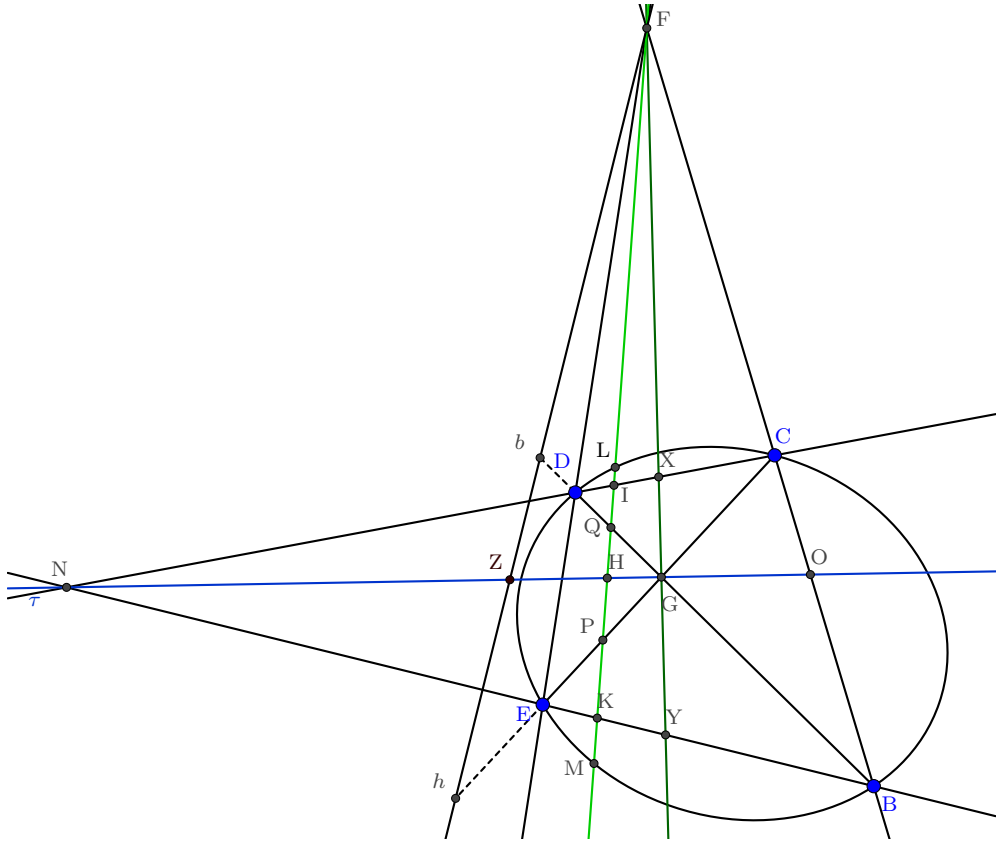
\begin{figure}[!ht]
\centering
\definecolor{ttqqqq}{rgb}{0.2,0.,0.}
\definecolor{qqttcc}{rgb}{0.,0.2,0.8}
\definecolor{qqccqq}{rgb}{0.,0.8,0.}
\definecolor{qqwuqq}{rgb}{0.,0.39215686274509803,0.}
\definecolor{uuuuuu}{rgb}{0.26666666666666666,0.26666666666666666,0.26666666666666666}
\definecolor{qqqqff}{rgb}{0.,0.,1.}
\begin{tikzpicture}[line cap=round,line join=round,x=0.7946210268948652cm,y=0.7946210268948652cm]
\clip(-5.28,-5.04) rectangle (11.08,8.88);
\draw [rotate around={161.3411074719956:(6.718011949781292,-1.614238588797848)},line width=0.8pt] (6.718011949781292,-1.614238588797848) ellipse (2.874329867174444cm and 2.4252535431644815cm);
\draw [line width=0.8pt,domain=-5.28:11.08] plot(\x,{(--0.1124-0.62*\x)/-3.3});
\draw [line width=0.8pt,domain=-5.28:11.08] plot(\x,{(--10.4752--1.36*\x)/-5.48});
\draw [line width=0.8pt,domain=-5.28:11.08] plot(\x,{(--14.1852-3.54*\x)/-0.54});
\draw [line width=0.8pt,domain=-5.28:11.08] plot(\x,{(-43.1888--5.52*\x)/-1.64});
\draw [line width=0.8pt] (4.12,0.74)-- (9.06,-4.16);
\draw [line width=0.8pt] (3.58,-2.8)-- (7.42,1.36);
\draw [line width=0.8pt,color=qqwuqq,domain=-5.28:11.08] plot(\x,{(--50.626487701184566-9.161210824141644*\x)/0.24191301396863363});
\draw [line width=0.8pt,color=qqccqq,domain=-5.28:11.08] plot(\x,{(--34.42614957661839-7.328806337066371*\x)/-0.5220213056541922});
\draw [line width=0.8pt,color=qqttcc,domain=-5.28:11.08] plot(\x,{(-7.568587535227991--0.1705862809068719*\x)/9.849530342753594});
\draw [line width=0.8pt,domain=-5.28:11.08] plot(\x,{(--29.559679923718694-9.20466226904343*\x)/-2.266942462107246});
\draw [line width=0.8pt,dash pattern=on 2pt off 2pt] (3.5362362794639584,1.3190368887908106)-- (4.12,0.74);
\draw [line width=0.8pt,dash pattern=on 2pt off 2pt] (3.58,-2.8)-- (2.136715498939216,-4.363558209482515);
\begin{scriptsize}
\draw [fill=qqqqff] (3.58,-2.8) circle (2.5pt);
\draw[color=qqqqff] (3.26,-2.95) node {$E$};
\draw [fill=qqqqff] (4.12,0.74) circle (2.5pt);
\draw[color=qqqqff] (3.9,1.27) node {$D$};
\draw [fill=qqqqff] (7.42,1.36) circle (2.5pt);
\draw[color=qqqqff] (7.56,1.69) node {$C$};
\draw [fill=qqqqff] (9.06,-4.16) circle (2.5pt);
\draw[color=qqqqff] (9.32,-4.45) node {$B$};
\draw [fill=black] (4.78,1.16) circle (1.0pt);
\draw[color=black] (4.58,1.45) node {$L$};
\draw [fill=uuuuuu] (-4.305596023130768,-0.8429907679821443) circle (1.5pt);
\draw[color=uuuuuu] (-4.42,-0.53) node {$N$};
\draw [fill=uuuuuu] (5.3020213056541925,8.488806337066372) circle (1.5pt);
\draw[color=uuuuuu] (5.58,8.65) node {$F$};
\draw [fill=uuuuuu] (5.543934319622826,-0.6724044870752725) circle (1.5pt);
\draw[color=uuuuuu] (5.68,-1.11) node {$G$};
\draw [fill=uuuuuu] (5.499792660032231,0.9992337724909043) circle (1.5pt);
\draw[color=uuuuuu] (5.64,1.29) node {$X$};
\draw [fill=uuuuuu] (5.613441988227601,-3.304649836494441) circle (1.5pt);
\draw[color=uuuuuu] (5.76,-3.01) node {$Y$};
\draw [fill=uuuuuu] (4.758630329873502,0.8599850922792641) circle (1.5pt);
\draw[color=uuuuuu] (4.94,0.71) node {$I$};
\draw [fill=uuuuuu] (4.4819897233532515,-3.0238514641898577) circle (1.5pt);
\draw[color=uuuuuu] (4.62,-2.73) node {$K$};
\draw [fill=uuuuuu] (4.78,1.16) circle (1.5pt);
\draw [fill=uuuuuu] (4.427994876718701,-3.781900551263631) circle (1.5pt);
\draw[color=uuuuuu] (4.1,-3.85) node {$M$};
\draw[color=qqttcc] (-4.83,-1.08) node {$\tau$};
\draw [fill=uuuuuu] (4.648375436240884,-0.6879148773104405) circle (1.5pt);
\draw[color=uuuuuu] (4.8,-0.41) node {$H$};
\draw [fill=uuuuuu] (8.011135168102461,-0.629674468247313) circle (1.5pt);
\draw[color=uuuuuu] (8.16,-0.33) node {$O$};
\draw [fill=uuuuuu] (4.708504931176668,0.15626029093812388) circle (1.5pt);
\draw[color=uuuuuu] (4.4,0.15) node {$Q$};
\draw [fill=uuuuuu] (4.574688966809071,-1.7224202859568385) circle (1.5pt);
\draw[color=uuuuuu] (4.32,-1.49) node {$P$};
\draw [fill=ttqqqq] (3.0350788435469465,-0.7158559319770581) circle (1.5pt);
\draw[color=ttqqqq] (2.8,-0.41) node {$Z$};
\draw [fill=uuuuuu] (3.5362362794639584,1.3190368887908106) circle (1.5pt);
\draw[color=uuuuuu] (3.2,1.59) node {$b$};
\draw [fill=uuuuuu] (2.136715498939216,-4.363558209482515) circle (1.5pt);
\draw[color=uuuuuu] (1.86,-4.07) node {$h$};
\end{scriptsize}
\end{tikzpicture}
\caption{L'involution de polarité $F,Z;b,h$ induite sur l'ordinale $(FZ)$.}\label{Traversale-09}
\end{figure}

Menons, de $F$, une tangente à la conique $\cC$. La traversale $\tau$ de $F$ passe alors nécessairement par le point $S$ où cette tangente touche la conique. L'involution de Desargues, toujours hyperbolique, admet donc $F$ et $S$ comme n{\oe}uds moyens doubles. Mais l'involution de polarité devient dégénérée, puisque ses deux n{\oe}uds moyens doubles se confondent en $S$. En quelque sorte, tous les points, sauf $S$, sont accouplés à $F$. L'arbre correspondant~«~est de l'espèce mitoyenne, dont l'entendement ne peut comprendre comment sont les proprietez que le raisonnement luy en fait conclure.~» Le mot mitoyen est soigneusement choisi par Desargues. En effet, si l'on considère une droite passant par $F$ et \textit{ne coupant pas} la conique $\cC$, alors l'involution de polarité est une involution elliptique, c'est-à-dire sans n{\oe}ud moyen double. Lorsque l'on fait varier cette droite et qu'elle vient couper la conique, l'involution devient hyperbolique. Entre ces deux situations, la droite est devenue à un moment donné tangente à la conique et, dans cette situation intermédiaire ou mitoyenne, l'involution dégénère. La figure \ref{Traversale-09} permet de se faire une idée de la chose. 

L'usage, dans le cas de la tangente, du mot d'arbre plutôt que de celui d'involution montre bien que Desargues, prudent mais clairvoyant, a compris que, sur une tangente, il n'y a pas à proprement parler d'involution de polarité\footnote{Rappelons que dans \cite{lenger} comme dans \cite{taton} il est écrit que Desargues découvre ici le troisième cas d'involution, dit parabolique, alors que d'une part Desargues ne lui attribue pas le nom d'involution et que, d'autre part, il n'existe pas d'involution parabolique.}. 

\begin{figure}[!ht]
\centering
\includegraphics[width=11cm]{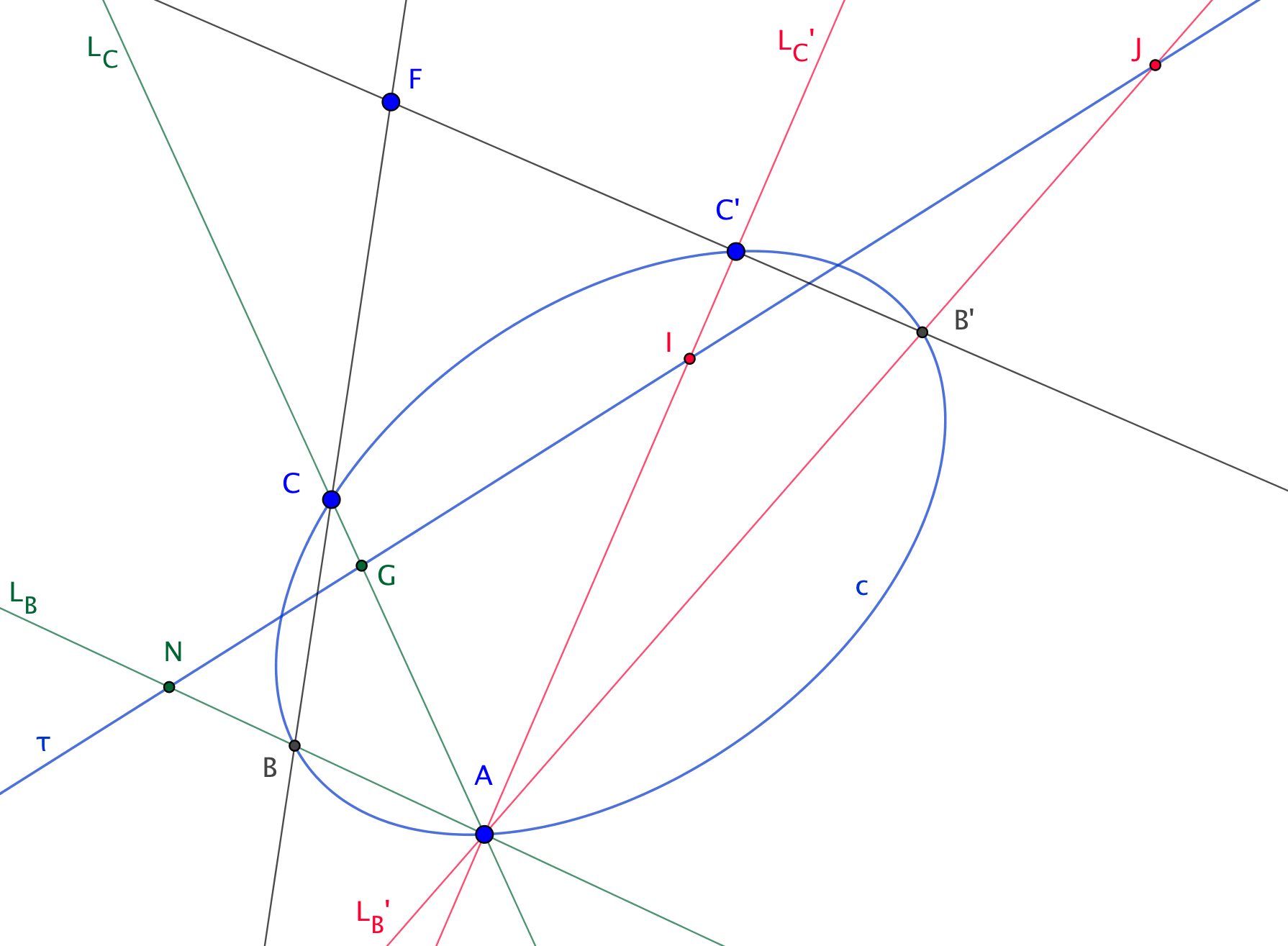}
\caption{Construction de l'involution de polarité sur $\tau$.}\label{Involution-Polaire-01}
\end{figure}

La profondeur de ses vues est confirmée par les lignes 22 à 25 de la page 22. Voici ce qu'il y dit. Soient $\cC$ une coupe de rouleau et $\tau$ une droite, $F$ le but de  l'ordonnance dont $\tau$ est la traversale eu égard à $\cC$. Soit $A$ un point sur la conique et $L_C,L_B$ deux droites de même ordonnance en $A$ et qui coupent $\cC$ en deux points $B$ et $C$ alignés avec $F$, c'est à dire situés sur une même droite de l'ordonnance en $F$. Alors ces deux droites coupent $\tau$ en deux points $G$ et $N$. Si l'on répète cette construction avec autant de couples de droites passant par $A$ et ayant la même propriété d'alignement avec $F$ de leurs points d'intersection avec $\cC$, nous obtenons autant de couples de points $N,G; I,J;\ldots$ d'une involution, qui n'est autre que l'involution de polarité.  Nous renvoyons à la figure \ref{Involution-Polaire-01} pour une illustration de cette construction qui, soit dit au passage, est très proche de celle utilisée pour construire des paramétrages rationnels d'une conique par le pinceau des droites passant par un point pris sur icelle. 

Ces lignes nous montrent que Desargues a bien en tête qu'une involution peut comprendre \emph{autant de couples de points} qu'il est nécessaire d'en considérer. En extrapolant un peu, on pourrait dire qu'il a entre les mains un \emph{nouvel objet mathématique,} à savoir une \emph{relation.} On peut penser que la difficulté à conceptualiser cette nouvelle notion fait partie des raisons qui ont pu l'empêcher de donner une présentation définitive de sa théorie de la traversale et de la polarité.

\subsection{La conclusion de Desargues~:~diamètres conjugués, asymptotes}
Nous allons pour finir rapporter intégralement les lignes que Desargues écrit en conclusion de la présentation de sa théorie de la traversale\footnote{p. 22, l. 46 à 55.}. Il revient alors à la géométrie dans l'espace et donc, implicitement, à ses préoccupations perspectivistes~:~

«~Mais voicy dans une proposition comme un assemblage abregé de tout ce qui precede.

Estant donnée de grandeur \& de position une quelconque coupe de rouleau à bord courbe E, D, C, B, pour assiette ou base d'un quelconque rouleau, dont le sommet soit aussi donné de position, \& qu'un autre plan en quelconque position aussi donnée coupe ce rouleau, \& que l'essieu 4,5, de l'ordonnance de ce plan de coupe avec le plan d'assiette ou base soit aussi donné de position, la figure qui vient de cette construction en ce plan de coupe est donnée d'espece \& de position, chacune de ses diametrales avec leurs distinctions de coniuguées \& essieux, comme encore chacune des especes de leurs ordonnées \& des touchantes à la figure, \& la nature de chacune; leurs ordonnances, avec les distinctions possibles, sont donnez tous de generation \& de position.~»

\begin{figure}[!ht]
\centering
\includegraphics[width=11cm]{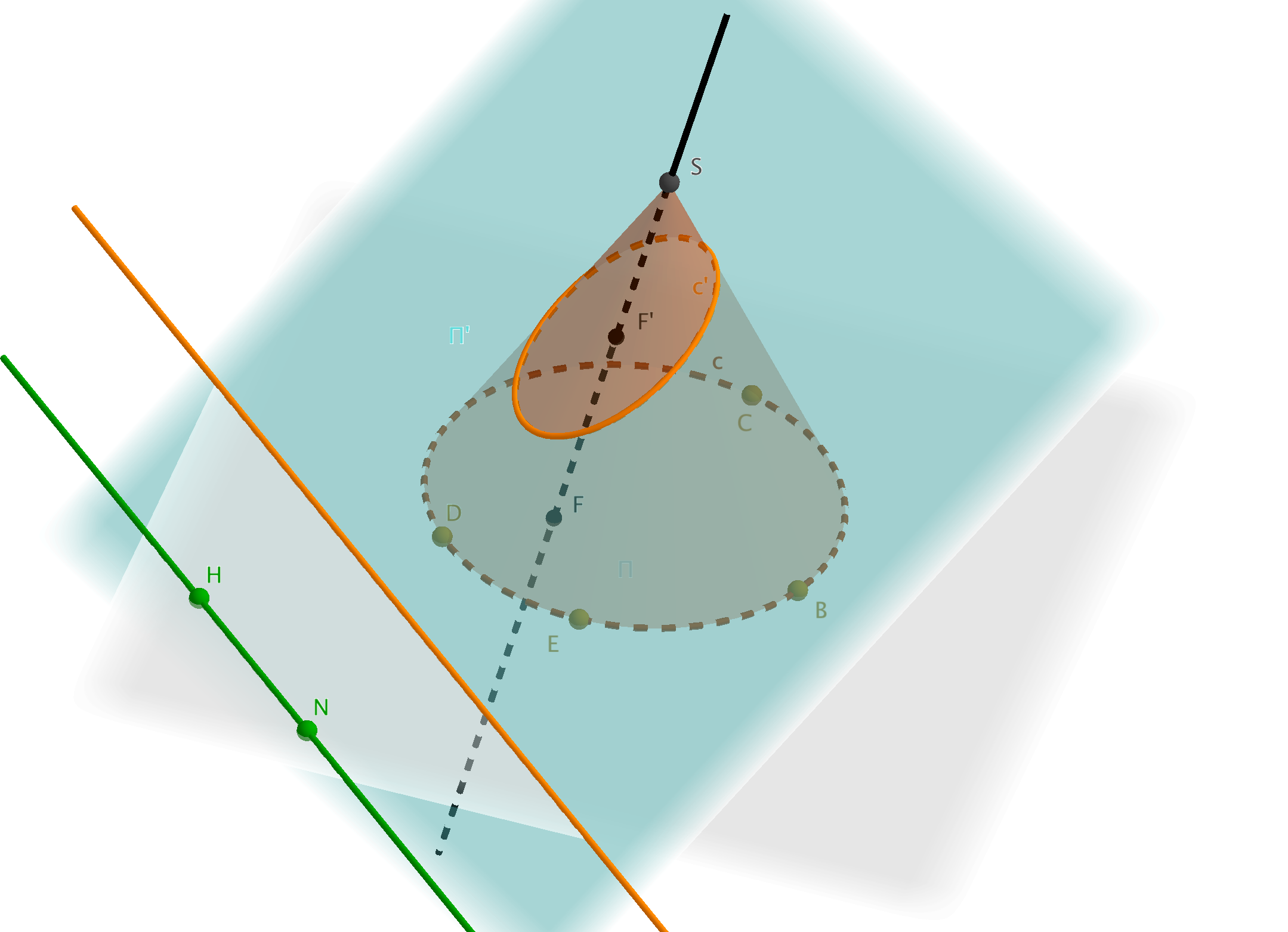}
\caption{Construction du centre d'une conique par perspective depuis une autre, grâce à la polaire.}\label{Conclusion-3D-01}
\end{figure}

On peut reformuler cela en disant que la polarité étant affaire d'incidence, comme l'a montré Desargues dans les lignes qui précèdent, elle se transfère d'un plan de coupe à un autre par simple projection centrale depuis le sommet du cône. En voici une illustration donnée par l'auteur juste après les lignes citées ci-dessus, et que nous reformulerons ainsi : donnons-nous un cône de sommet $S$ et un plan $\Pi$ coupant ce cône selon une conique $\cC$. Considérons un second plan $\Pi'$ coupant ce même cône suivant une conique $\cC'$ et traçons par $S$ le plan $\Pi"$ parallèle à $\Pi'$. Le plan $\Pi$ coupe les plans $\Pi'$ et $\Pi"$ selon les droites $45$ et $NH$ qui, par construction, sont parallèles. Si l'on s'imagine représenter en perspective le plan $\Pi'$ dans le plan $\Pi$ depuis le point $S$, alors la droite $NH$ est la droite à l'infini de ce plan $\Pi'$ dans cette représentation. De la sorte, si l'on considère le point $F$ dont la droite $NH$ est la traversale eu égard à la conique $\cC$, et si l'on trace la droite $FS$ qui va donner en le plan $\Pi'$ le point $F'$, ce dernier point est le \emph{centre} de la conique $\cC'$. Ainsi la droite $FS$ est-elle l'essieu d'une ordonnance de plans qui coupent le plan $\Pi'$ selon des diamètres, ou \emph{diamétrales,} de la conique $\cC'$. 

Mais il y a mieux et Desargues poursuit en montrant comment on peut dans sa théorie retrouver les diamètres conjugués. Choisissons $Z$ sur la droite $FH$, considérons sa traversale $\tau$ et supposons, quitte à bouger le point $H$ sur $FH$, que $H$ est le point d'intersection de $FH$ et de $\tau$. De la sorte, $Z,H$ est un couple de points dans l'involution de polarité induite par $\cC$ sur $FH$. Traçons les deux plans $SFH$ et $SFZ$~:~ils coupent $\Pi'$ en deux diamètres qui sont \textit{conjugués\footnote{p. 23, l. 1 à 5.}.} Desargues va plus loin et considère maintenant les deux droites $SZ$ et $SH$ comme des essieux d'ordonnances de plans qui engendrent deux pinceaux de droites dans le plan $\Pi'$, formant des ordonnées conjuguées entre elles (car provenant de diamètres conjugués). 

Ces considérations sont valides que ces ordonnées soient à distance finie ou infinie, et qu'elles rencontrent, ne rencontrent pas ou simplement touchent la conique $\cC'$. Si en outre $\cC'$ est une hyperbole, alors les deux asymptotes sont du corps des droites de l'ordonnance au centre. Elles sont conjuguées, et peuvent donc être vues à la fois comme tangentes en les points à distance infinie de l'hyperbole et comme diamètres conjugués\footnote{p. 23, l. 10 à 15.}. Comme l'écrit Taton en note de la page 159 de \cite{taton}~:~«~Cette théorie projective des diamètres conjugués et des asymptotes, celles-ci étant considérées comme diamètres conjugués et tangentes à la conique à l'infini, est une des parties les plus élégantes du traité de Desargues.~»

Nous finirons cette section par quelques remarques ; aux lignes 26 à 31 de la page 23, Desargues introduit pour la première fois (c'est-à-dire, avant de le faire rétrospectivement dans les \textit{Advertissements}) la notion d'ordinale, qu'il réserve aux tangentes vues comme du corps des ordonnées d'un diamètre. Ce n'est donc qu'après les développements décrits ci-dessus concernant l'hyperbole et ses asymptotes qu'il réalisera que cette notion est valide en toute généralité, expliquant ses ajouts ultérieurs mentionnés dans les \textit{Advertissements,} justifiant une fois de plus l'utilité de la lecture de \textit{l'original} du \textit{Brouillon Project} en tenant compte de l'articulation entre le corps du texte et ses corrections ultérieures. Enfin, la \textit{Figure 13} du manuscrit de la Hire\footnote{Reproduite ici à la figure \ref{LaHire-01}}, à laquelle ce passage renvoie,  qui est celle qu'il a déjà utilisée pour illustrer le théorème d'involution, montre la droite $45$, parallèle à $NH$. On est en droit de penser qu'il y a dans ces lignes du \textit{Brouillon Project} la source de l'inspiration de la Hire pour ses \textit{Planiconiques,} où il montre comment générer une conique à partir d'un cercle au moyen d'une homologie. Nous renvoyons à la fin de l'article \cite{cracewski} pour une analyse plus poussée allant en ce sens.
\section{Conclusion}
En 1822, dans l'introduction de son \textit{Traité des propriétés projectives des figures\footnote{{\it Voir} \cite{poncelet}. Nous citons l'édition de 1822 pour la pagination.},} Jean-Victor Poncelet décrit ce qui peut passer pour un programme de refondation de la géométrie. Il souhaite donner à la géométrie synthétique héritée des anciens la même capacité à raisonner en toute généralité que possède la géométrie analytique grâce à son usage de l'algèbre. Ainsi écrit-il, aux pages XI et XII~:~

«~L'Algèbre emploie des signes abtraits, elle représente des grandeurs absolues par des caractères qui n'ont aucune valeur par eux-mêmes, et qui laissent à ces grandeurs toute l'indétermination possible ; par suite, elle opère et raisonne forcément sur les signes de non-existence comme sur des quantités toujours absolues, toujours réelles [\ldots]. Le résultat doit donc lui-même participer de cette généralité, et s'étendre à tous les cas possibles, à toutes les valeurs des lettres qui y entrent ; de là aussi ces formes extraordinaires, ces êtres de raison, qui semblent l'apanage exclusif de l'Algèbre. 

Or on est conduit à toutes ces conséquences, non-seulement quand on emploie les signes et les notations de l'Algèbre, mais aussi toutes les fois qu'en raisonnant sur des grandeurs quelconques on fait abstraction de leurs valeurs numériques et absolues ;  en un mot, toutes les fois qu'on emploie le raisonnement sur des grandeurs \textit{indéterminées,} c'est-à-dire le \textit{raisonnement} purement \textit{implicite.}~» 

Il poursuit en évoquant le \textit{principe de continuité,} qui l'amène à considérer que deux droites parallèles se coupent, en un point, et un seul, situé à une distance «~plus grande que toute distance donnée.~» On reconnait là ce qu'écrivait Desargues au début du \textit{Brouillon Project,} presque deux siècles auparavant~:~

«~Pour donner à entendre l'espece de position d'entre plusieurs droictes en laquelle elles sont toutes paralelles entre elles, il est icy dict, que toutes ces droictes sont entre elles d'une mesme ordonnance, dont le but est à distance infinie en chacune d'elles d'une part \& d'autre. [\ldots] Ainsi deux quelconques droictes en un mesme plan sont entre elles d'une mesme ordonnance, dont le but est à distance ou finie, ou infinie.\footnote{p. 1, l. 22 à 24 et 28 à 29.}~» 

Poncelet évoque aussi, comme conséquence du principe de continuité, la possiblité d'étendre un énoncé à des cas où certains points se confondent ou s'éloignent à l'infini, rappelant les termes de Desargues~:~

«~Chacun pensera ce qui lui semblera convenable ou de ce qui est icy deduit, ou de la maniere de le deduire, \& verra que la raison essaye à cognoistre des quantitez infinies d'une part: ensemble de celles qui s'apetissent iusques à reduire leurs deux extremitez opposées en une seule, \& que l'entendement s'y pert, non seulement à cause de leurs inimaginables grandeur et petitesse, mais encore à cause que le raisonnement ordinaire le conduit à en conclure des proprietez, d'où il est incapable de comprendre comment, c'est qu'elles sont.\footnote{p. 1, l. 5 à 10.}~»

Toujours à la recherche de la plus grande généralité évoquée un peu avant, Poncelet poursuit son introduction en montrant que, pour atteindre son objectif, il faut se débarasser de~«~l'usage trop fréquent et trop étendu du mécanisme des proportions, qui n'est au fond qu'un calcul déguisé,~» pour se concentrer sur des propriétés qui~«~auront toute la généralité, toute l'indétermination possible.~» Ces propriétés comprennent celles qui se conservent en prenant diverses \textit{transversales} à une figure donnée. Ainsi,~«~la Théorie des transversales [$\cdots$] n'est à proprement parler que la Théorie des coordonnées prises sur une même droite, et réduite par conséquent à un plus grand degré de simplicité.\footnote{\cite{poncelet}, p. XIX.}~»

Il conclut en écrivant que cette théorie des transversales n'est qu'un corollaire des «~principes de la méthode des projections,~» autrement dit, de la méthode qui se base exclusivement sur des propriétés qui sont vraies~«~à la fois pour une figure donnée et pour toutes ses projections ou perspectives.~» Ainsi ces propriétés doivent être distinguées~«~de toutes les autres par le nom générique de \textit{propriétés projectives,} qui en rappellent, de manière abrégée, la véritable nature.~» Outre le même mot de «~nature~» employé par Poncelet et Desargues ({\it voir} plus haut), ces lignes font résonner en nous l'écho du \textit{Brouillon Project} et de sa disposition de points en involution sur une transversale, disposition qui subsiste par projection centrale, ou ramée, et qui va révolutionner l'étude des coniques grâce à des \textit{traversales} qui ne diffèrent que par une lettre dans leur nom des \textit{transversales.}

Un peu plus loin dans son introduction, après quelques considérations historiques sur la géométrie antique, Poncelet écrit ceci~:~«~Desargues, ami de l'illustre Descartes, et dont celui-ci faisait le plus grand cas comme géomètre ;  Desargues, qu'on peut appeler, à plus d'un titre, le Monge de son siècle, que les biographes n'ont point assez connu, ni assez compris ;  Desargues, enfin, que des contemporains, indignes du beau titre de géomètre, ont noirci, persécuté, dégouté, pour n'avoir pu se mettre à la hauteur de ses idées et de son génie, fut, je crois, le premier, d'entre les modernes qui envisagea la Géométrie sous le point de vue général que je viens de faire connaitre.\footnote{p. XXV et XXVI.}~» 

On ne peut qu'approuver Poncelet, dont la conviction était basée sur des sources très précaires\footnote{Principalement les \textit{Advis charitables} de Jean de Beaugrand, quelques lettres de Descartes ou de Mersenne, et le traité de perspective d'Abraham Bosse.}, lorsqu'il déclare en quelque sorte Desargues «~père de la géométrie projective~». L'usage indifférencié de points à distance finie ou infinie, l'utilisation de transversales à des figures pour leur étude et l'introduction d'une notion invariante par projection sont en effet la marque de tous les traités de géométrie projective qui verrons le jour après 1822. La lecture d'autres passages de l'ouvrage de Poncelet, comme de la section II du chapitre I par exemple, donne d'ailleurs à celui qui s'est penché sérieusement sur le \textit{Brouillon Project} de Desargues un curieux sentiment de familiarité. 

Notre analyse détaillée de l'original du \textit{Brouillon Project} nous fait conclure en outre que Desargues est également l'inventeur de la théorie de la polarité. Qu'il ait vu là une nouvelle \textit{théorie} ne fait pour nous pas de doute. Cependant les difficultés qu'il a pu rencontrer à sa rédaction, voire l'attitude hostile de ses contemporains, l'ont empêché d'en donner une exposition claire et convaincante. Il est à ce titre intéressant de comparer le \textit{Brouillon} avec le petit \textit{Mémoire sur les lignes du second ordre\footnote{\textit{voir} \cite{brianchon}.}} de Charles-Julien Brianchon qui, en 1817, redécouvre la théorie de la polarité, introduisant au passage les mots de «~pôle~» et de «~polaire~». Il y utilise des méthodes très proches de celles de Desargues, basées sur un usage intensif du théorème de Ménélaüs, et y reprend les idées générales de Carnot sur l'importance des \textit{transversales,} droites quelconques coupant une figure que l'on souhaite étudier\footnote{Nous renvoyons au texte de l'habilitation à diriger des recherche \cite{nabonnand-habil} de Philippe Nabonnand pour une étude détaillée de la naissance de la polarité au début du XIX$^e$ siècle.}. Brianchon surmonte la difficulté rencontrée par Desargues en introduisant la notion de pôle d'une droite, oubliant l'ordonnance des droites en ce point, en faisant ainsi clairement apparaître une dualité, ou correspondance, entre des points et des droites. Nous renvoyons le lecteur à la page 13 de ce traité, où apparaissent les conditions de l'involution, à ses pages 14 et 15, où est énoncée une version du théorème d'involution de Desargues, suivi d'une remarque sur l'inutilité de développer outre mesure la théorie des diamètres, celle-ci étant contenue dans celle des polaires. Mais, d'une manière générale, on peut considérer ce texte comme moins ambitieux que ne l'était le \textit{Brouillon Project.} 

Enfin, outre nous montrer l'indéniable génie de Desargues, la lecture de l'original du \textit{Brouillon Project} nous donne l'occasion unique d'entrer dans \textit{l'atelier du géomètre,} ce qui est extrêmement rare compte tenu de la propension des mathématiciens à ne livrer à l'impression que des textes achevés dont ont été exclues les tentatives infructueuses. La délicate articulation entre le corps du texte et les notes rajoutées \textit{a posteriori} permet presque de suivre pas à pas le développement de la pensée de Desargues. Cela justifie à notre avis la nécessité de l'établissement d'une nouvelle édition critique de son traité.


\section*{Appendice~:~la théorie de la polarité dans le langage contemporain de la géométrie projective}
Les traités modernes de géométrie projective ou de géométrie des coniques définissent directement les coniques à partir des formes quadratiques, approche qui donne immédiatement la notion de polaire, permettant d'arriver très rapidement et simplement aux résultats les plus profonds concernant les coniques, à l'exception des plus difficiles démontrés au dix-neuvième siècle comme le «~porisme~» de Poncelet\footnote{La démonstration de ce théorème fait aujourd'hui appel à la théorie des cubiques, plus précisément des courbes elliptiques.} ou le théorème d'énumération de Chasles\footnote{Chasles y dénombre les coniques tangentes à cinq coniques données~:~il y en a au plus 3264.}.

La théorie des coniques ne faisant plus partie du cursus mathématique depuis quelques décennies, nous allons ici, pour la commodité du lecteur, en rappeler quelques-uns des grands traits, en nous centrant sur la théorie de la polarité. Cet appendice se veut aussi un plaidoyer pour le retour de la théorie projective des coniques dans les programmes de mathématiques des premiers et seconds cycles universitaires, comme mise en pratique de la théorie des formes bilinéaires et de la dualité vectorielle.

Mais commençons par une mise en bouche~:~selon la définition antique, une surface conique ou cône est l'union des droites passant par un sommet $S$ et par les points d'un cercle inclus dans un plan ne contenant pas $S$. Choisissons un repère affine dont $S$ soit l'origine et dont les axes de coordonnées, que l'on peut supposer orthonormés, sont tels que le cercle soit décrit par les deux équations $z=1$ et $x^2+y^2=r^2$. L'équation du cône s'obtient alors simplement en homogénéïsant l'équation du cercle~:~$x^2+y^2-r^2z^2=0$. À la vue de cette équation, il semble naturel d'introduire l'étude des coniques par l'étude des formes quadratiques.
\subsection{Formes quadratiques, formes bilinéaires, dualités}
Dans toute la suite et afin de nous placer dans le cadre classique de la théorie des coniques, nous considèrerons que $E$ est un espace vectoriel réel de dimension 3, même si une bonne partie de notre propos s'étend sans peine à bien d'autres cas. Nous noterons $E^*$ le dual de $E$.

On dit qu'une application $q$ de $E$ dans $\RR$ est une \emph{forme quadratique} si, exprimée en coordonnées dans une base (et cela sera alors vrai dans toute base), $q$ est un polynôme homogène de degré $2$. Notons $\cQ(E)$ l'espace vectoriel des formes quadratiques sur $E$. 

Une application $\phi$ de $E\times E$ dans $\RR$ est une \emph{forme bilinéaire} si pour tout vecteur $\vec{v}$ de $E$, les deux applications
\[
\vec{w}\mapsto \phi(\vec{v},\vec{w})\;\mbox{et}\;\vec{w}\mapsto \phi(\vec{w},\vec{v})
\]
sont des formes linéaires.  Notons $\Bilin(E)$ l'espace vectoriel des formes linéaires sur $E$.

Soit $t$ un élément du produit tensoriel $E^*\otimes_{\RR}E^*$ ; on peut écrire $t$ comme combinaison linéaire de tenseurs purs, à savoir
\[
t=\sum_i f_i\otimes g_i, \; f_i,g_i\in E^*,
\]
et il induit donc une forme bilinéaire $\phi_t$ sur $E$ de la manière suivante~:~
\[
\phi_t(\vec{v},\vec{w})=\sum_i f_i(\vec{v})g_i(\vec{w}).
\]
On réalise ainsi un isomorphisme naturel entre $\Bilin(E)$ et $E^*\otimes E^*$\footnote{C'est une tautologie si l'on se base sur la propriété universelle du produit tensoriel.}.

Une forme bilinéaire $\phi\in\Bilin(E)$ est \emph{symétrique} si pour tous vecteurs $\vec{v},\vec{w}$ de $E$ on a
\[
\phi(\vec{v},\vec{w})=\phi(\vec{w},\vec{v}).
\]
Elle est \emph{alternée\footnote{Rappelons que comme nous travaillons sur le corps des nombres réels, qui est de caractéristique différente de $2$, le caractère alterné est équivalent au caractère antisymétrique.}} si pour tout vecteur $\vec{v}$ on a 
\[
\phi(\vec{v},\vec{v})=0.
\]
Les formes bilinéaires symétriques forment un sous-espace noté $\Bilin^s(E)$ et les formes alternées un sous-espace noté $\Bilin^a(E)$. Toute forme bilinéaire sur $E$ se décompose de manière unique en la somme d'une forme symétrique et d'une forme alternée, de sorte que
\[
\Bilin(E)=\Bilin^s(E)\oplus \Bilin^a(E),
\]
et cette décomposition se traduit dans l'espace $E^*\otimes E^*$ par une décomposition analogue~:~
\[
E^*\otimes E^*=(E^*\odot E^*)\oplus (E^*\vect E^*).
\]

Si $\phi$ est une forme bilinéaire sur $E$, elle définit sur $E$ une forme quadratique $q_\phi$ en posant
\[
q_\phi(\vec{v})=\phi(\vec{v},\vec{v}).
\]
L'application $\phi\mapsto q_\phi$ est linéaire et son noyau est l'espace des formes alternées. Ainsi elle réalise un isomorphisme naturel entre $\Bilin^s(E)$ et $\cQ(E)$. On peut expliciter sa réciproque~:~si $q\in \cQ(E)$, la \emph{forme polaire} de la forme quadratique $q$ est la forme bilinéaire symétrique $\phi_q$ définie par
\[
\phi_q(\vec{v},\vec{w})=\frac{1}{2}\left[q(\vec{v}+\vec{w})-q(\vec{v})-q(\vec{w})\right].
\]

Rappelons que si $f\in E^*$ et $\vec{v}\in E$, il est d'usage de noter $f(\vec{v})$ sous la forme $\la f | \vec{v}\ra$, définissant ainsi les \emph{crochets de dualité} chers aux amateurs de mécanique quantique. Nous appellerons \emph{dualité sur $E$} toute application linéaire $d$ de $E$ dans $E^*$ qui est telle que, pour tous vecteurs $\vec{v},\vec{w}$, on ait
\[
\la d(\vec{v}) | \vec{w}\ra =\la d(\vec{w}) | \vec{v} \ra.
\]
Une dualité définit ainsi une forme bilinéaire symétrique $\phi_d$ et, réciproquement, toute forme bilinéaire symétrique $\phi$ définit une dualité $d_\phi$ par
\[
\la d_\phi(\vec{v}) | \vec{w} \ra = \phi(\vec{v},\vec{w}).
\]
Notons $\cD(E)$ l'espace vectoriel des dualités sur $E$, qui est donc naturellement isomorphe à l'espace $\Bilin^s(E)$, ce que nous aurions pu déduire directement de l'isomorphisme canonique entre $\Hom(E,E^*)$ et $E^*\otimes E^*$.

Nous pouvons résumer ce qui précède de la manière suivante~:~les espaces $\cQ(E)$ (des formes quadratiques sur $E$), $\Bilin^s(E)$ (des formes bilinéaires symétriques), $E^*\odot E^*$ (des tenseurs symétriques) et $\cD(E)$ (des dualités) sont isomorphes les uns aux autres.

Rappelons que le \textit{rang} d'une forme bilinéaire $\phi$ est, par définition, celui de la dualité associée $d_\phi$. Si $\phi$ est de rang maximal, ici $3$, on dit qu'elle est \textit{non-dégénérée.} La dualité associée est alors un isomorphisme. Comme le bidual $E^{**}$ est canoniquement isomorphe à $E$, toute forme bilinéaire symétrique sur $E$ en induit une sur $E^*$, ce qui permettra par exemple de construire les \emph{coniques tangentielles.}
\subsection{Cône isotrope et conique}
La théorie classique des coniques s'intéresse aux sections d'un cône. Avant de nous pencher sur le mot «~section~», voyons comment une forme quadratique permet de définir, de manière naturelle, un «~cône~». Si $q$ est une forme quadratique sur $E$, son \emph{cône isotrope} est l'ensemble des vecteurs $\vec{v}$ de $E$ tels que $q(\vec{v})=0$. De l'homogénéïté de l'application $q$, se traduisant par le fait que si $\lambda\in \RR$ alors $q(\lambda \vec{v})=\lambda^2q(\vec{v})$, on peut déduire deux faits importants concernant ce cône isotrope. Tout d'abord, il s'agit bien d'un cône~:~s'il n'est pas réduit au vecteur nul, il est une union de droites vectorielles de $E$. Ensuite, si $\mu$ est un nombre réel non-nul quelconque, le cône isotrope de $\mu q$ est le même que celui de $q$. Dit autrement, le cône isotrope de $q$ ne dépend que du point $\bar{q}$ que définit $q$ dans l'espace projectif $\PP(\cQ)$ de $\cQ$. Ce point $\bar{q}$ s'appelle la \emph{classe conforme} de $q$. Nous noterons donc le cône isotrope de $q$ sous la forme $\CI(\bar{q})$. 

Pour former une section conique avec le cône $\CI(\bar{q})$, nous pourrions prendre un plan affine dans $E$ et considérer l'intersection de $\CI(\bar{q})$ avec celui-ci. Afin d'unifier toutes les sections possibles obtenues de cette manière, il est naturel de considérer l'image du cône isotrope dans l'espace projectif $\PP(E)$. Nous noterons cette image $\cC_{\bar{q}}$ et l'appellerons la \emph{conique associée à la (classe de la) forme quadratique $q$.} De cette manière, réaliser une conique classique comme \emph{section} du \emph{cône} isotrope revient simplement à choisir une carte affine dans $\PP(E)$, c'est-à-dire à y choisir une droite que l'on considèrera comme la droite à l'infini du plan affine coupant $\CI(\bar{q})$. 

Si le cône isotrope de $q$ est réduit à $\{\vec{0}\}$ alors sa conique associée est vide. Une telle conique pourrait aussi être réduite à un point, à une droite ou à l'union de deux droites, si la forme quadratique $q$ est dégénérée. Supposons maintenant que  $q$ est non dégénérée. Selon le théorème d'inertie de Sylvester, il existe une \emph{base d'inertie} pour la forme quadratique $q$ \footnote{Cette terminologie est d'origine mécanique : une telle base permet d'obtenir les axes principaux d'inertie d'un solide dont on connaît la matrice d'inertie, qui est symétrique et non dégénérée}, c'est-à-dire une base $\cB_i=\{\vec{e}_0, \vec{e}_1, \vec{e}_2\}$ telle que, en coordonnées dans cette base, $q$ s'exprime de la manière suivante~:~
\[
q(\vec{v})=\epsilon_0 x^2+\epsilon_1 y^2 +\epsilon_2 z^3,
\]
où $(\epsilon_0,\epsilon_1,\epsilon_2)\in\{(1,1,1), (1,1,-1), (1,-1,-1), (-1,-1,-1)\}.$ La \emph{signature} de la forme quadratique vaut alors respectivement, dans ces différents cas~:~$(3,0), (2,1), (1,2)$ et $(0,3)$. Notons que comme l'on ne s'intéresse qu'au cône isotrope de $q$ et, plus généralement, à des objets s'obtenant par l'annulation de quantités dépendant linéairement de $q$, on peut restreindre notre étude aux cas où la signature vaut $(3,0)$ et $(2,1)$. Dans le premier cas, la forme quadratique $q$ est définie positive et donc $\cC_{\bar{q}}=\emptyset.$ C'est donc le second qui, géométriquement, donne naissance à une véritable conique, image dans l'espace projectif du cône usuel d'équation dans la base d'inertie $x^2+y^2-z^2=0$. Nous dirons alors que la conique $\cC_{\bar{q}}$ est une \emph{conique non dégénérée.} Comme dans le plan affine d'équation $z=1$ la trace du cône isotrope est le cercle d'équation $x^2+y^2=1$, il est facile de voir qu'une conique non dégénérée est homéomorphe à un cercle et qu'elle sépare $\PP(E)$ en deux composantes connexes. 

Remarquons que si la signature vaut $(1,1)$ ($q$ est donc dégénérée et de rang $2=1+1$), alors dans une base d'inertie la forme quadratique prend la forme suivante~:~
\[
q(\vec{v})=x^2-z^2=(x-z)(x+z),
\]
et ainsi le cône isotrope est réunion de deux plans, donnant une conique $\cC_{\bar{q}}$ qui est la réunion de deux droites projectives~:~c'est ce que nous appellerons dans la suite une \emph{conique dégénérée.} La droite commune à ces deux plans est constituée de vecteurs $\vec{v}$ qui non seulement sont isotropes au sens où $q(\vec{v})=\phi_q(\vec{v},\vec{v})=0$, mais qui sont en fait qualifiés de \emph{totalement isotropes~:~}pour tout vecteur $\vec{w}\in E$, on a $\phi_q(\vec{v}, \vec{w})=0$. Le point de $\PP(E)$ correspondant à cette droite est le \emph{point double} de la conique dégénérée $\cC_{\bar{q}}$, à savoir l'intersection des deux droites qui la constituent.

Notons pour finir que dans les deux cas que nous venons de traiter, à savoir dans les cas où la signature vaut $(2,1)$ ou $(1,1)$, on peut aisément démontrer que $\cC_{\bar{q'}}=\cC_{\bar{q}}$ si seulement si $\bar{q}=\bar{q'}$. Dit autrement, se donner une conique (dégénérée ou non dégénérée) est équivalent à se donner une classe conforme de forme quadratique sur $E$.
\subsection{Polarité}
Soit $q$ une forme quadratique sur $E$, que nous supposerons de signature égale à $(2,1)$ ou $(1,1)$, de sorte que la conique associée est soit non dégénérée et homéomorphe à un cercle, soit dégénérée et égale à l'union de deux droites se coupant en le point double de la conique. 

Soit $m\in \PP(E)$ un point qui ne soit pas un point double de $\cC_{\bar{q}}$. Ce point $m$ correspond à la droite vectorielle $D_m$ de $E$ qui, par hypothèse, ne contient que le vecteur nul comme vecteur totalement isotrope. Soit $\vec{m}$ un vecteur dirigeant la droite $D_m$ ; par analogie avec le produit scalaire usuel, on définit le $\bar{q}$-orthogonal de $D_m$ par
\[
D_m^\perp = \{\vec{w}\in E, \phi_q(\vec{m},\vec{w})=0\}.
\]
Il est évident que $D_m^\perp$ ne dépend que de $D_m$ et pas du vecteur choisi pour diriger $D_m$, et ne dépend en outre que de la classe conforme de $q$. De l'hypothèse sur $m$ découle que $D_m^\perp$ est un plan vectoriel de $E$, définissant une droite projective de $\PP(E)$ que nous noterons $m^{\perp}$~:~c'est la \emph{polaire du point $m$ eu égard à la conique $\cC_{\bar{q}}$.} La polaire de $m$ est bien définie relativement à la conique, puisqu'elle ne dépend que de la classe conforme qui nous donne la conique. En utilisant la dualité associée $d_\phi$, ainsi que l'identification naturelle entre l'espace $\PP(E)^*$ des droites projectives de $\PP(E)$ et l'espace  projectif $\PP(E^*)$ de $E^*$, il est immédiat (c'est-à-dire qu'il s'agit d'une tautologie) de voir que $m\mapsto m^\perp$ est une homographie de $\PP(E)$ dans $\PP(E)^*$.

Par converse, si l'on se donnne une droite projective $\delta$ de $\PP(E)$ de sorte que $\delta$ ne passe pas par l'éventuel point double de $\cC_{\bar{q}}$ et si l'on note $P_\delta$ le plan vectoriel de $E$ correspondant à $\delta$, on peut considérer le $q$-orthogonal de $P_\delta$, qui est une droite vectorielle donnant un point de $\PP(E)$ que l'on note $\delta^\perp$ et qui est appelé le \emph{pôle de la droite $\delta$ eu égard à la conique $\cC_{\bar{q}}$.} La correspondance $\delta \mapsto \delta^\perp$ est une homographie de $\PP(E)^*$ dans $\PP(E)$ et de la symétrie de la forme polaire $\phi_q$ on déduit qu'il s'agit de la correspondance réciproque de celle associant à un point sa polaire. Dit autrement nous avons les identités suivantes~:~
\[
(m^\perp)^\perp = m,\; (\delta^\perp)^\perp =\delta.
\]

Supposons maintenant que la conique soit non dégénérée, de sorte que dans une base d'inertie une équation du cône isotrope s'écrive
\[
x^2+y^2-z^2=0.
\]
La forme quadratique $q$, vue comme fonction, est différentiable en tout vecteur $\vec{v}$ et un calcul élémentaire montre que sa différentielle $Dq_{\vec{v}}$ au vecteur $\vec{v}$ satisfait à
\[
Dq_{\vec{v}}(\vec{w})=2\phi_q(\vec{v},\vec{w})
\]
pour tout vecteur $\vec{w}\in T_{\vec{v}}E=E$. Ainsi est-il clair que si $m$ est un point pris sur la conique $\cC_{\bar{q}}$, alors $m^\perp$ n'est jamais que la tangente à $\cC_{\bar{q}}$ en $m$. En fait $m\in m^\perp$ si et seulement si $m\in \cC_{\bar{q}}$ et réciproquement $\delta^\perp \in \delta$ si et seulement si $\delta$ est une tangente à la conique. 
\subsection{Polaire et involution}
Reprenons les hypothèses et notations de la section précédente. Si $\delta$ est une droite projective qui ne passe pas par l'éventuel point double de $\cC_{\bar{q}}$ et si l'on considère un point $m$ sur $\delta$, ce point $m$ admet lui-même une polaire $m^\perp$ qui, si $\delta$ n'est pas une tangente à la conique, est une droite différente de $\delta$, de sorte qu'elle coupe $\delta$ en un point $m'$. Un petit calcul dans une base d'inertie montre que $m\mapsto m'$ définit une homographie de $\delta$ dans elle-même. De la symétrie de $q_\phi$ on déduit que c'est une involution, que nous noterons $\pol_{\delta,\bar{q}}$ et nommerons \emph{involution de polarité induite sur $\delta$ par la conique $\cC_{\bar{q}}$.} Il est en effet évident qu'elle ne dépend que de la classe conforme de $q$, c'est-à-dire finalement de $\cC_{\bar{q}}$. Si $\delta$ est une tangente à la conique, alors cette involution n'est pas définie. Nous pourrions dire qu'elle est complètement dégénérée, envoyant tout point différent du point de tangence $m$ sur le point $m$ et le point $m$ sur la droite entière. Il faudrait pour donner un sens précis à cela procéder à un éclatement au point $m$, ce que nous ne développerons pas ici. 
\subsection{Coniques dans une carte affine}
Supposons ici que $q$ est non dégénérée et considérons une droite projective $\delta$. Trois cas peuvent se produire concernant les rencontres de cette droite et de la conique $\cC_{\bar{q}}$.
\begin{enumerate}
\item Si la droite $\delta$ ne rencontre pas $\cC_{\bar{q}}$, alors dans le plan affine $\PP(E)-\delta$ dont $\delta$ est la droite à l'infini, la conique $\cC_{\bar{q}}$ est une \emph{ellipse.} L'involution de polarité induite sur $\delta$ est sans point fixe et porte le nom d'\textit{involution elliptique.}
\item Si la droite $\delta$ rencontre la conique $\cC_{\bar{q}}$ en deux points distincts $a,a'$, alors dans le plan affine $\PP(E)-\delta$ la conique est une \emph{hyperbole} dont les deux asymptotes rencontrent la droite à l'infini $\delta$ en les points $a,a'$. Ces deux asymptotes sont par ailleurs les tangentes à $\cC_{\bar{q}}$ en les points $a,a'$ respectivement. L'involution de polarité induite sur $\delta$ a deux points fixes $a$ et $a'$ et porte le nom d'\textit{involution hyperbolique.}
\item Si la droite $\delta$ ne rencontre $\cC_{\bar{q}}$ qu'en un seul point, elle lui est tangente. Dans le plan affine $\PP(E)-\delta$ la conique est une parabole et l'involution induite sur $\delta$ n'en étant pas une, elle ne porte pas de nom.
\end{enumerate}
On retrouve ainsi l'origine de la terminologie concernant les homographies d'une droite projective suivant qu'une telle homographie admet deux points fixes (elle est hyperbolique), aucun point fixe (elle est elliptique), ou un seul (elle est parabolique). 

Nous pourrions poursuivre ainsi et retrouver les résultats et notions présentes dans le \textit{Brouillon,} comme par exemple la théorie des diamètres et des asymptotes, ou bien les constructions à la règle seule des pôles et des polaires, mais nous laissons au lecteur le plaisir de le faire découvrir à ses étudiants.

\vspace*{.5cm}
\noindent \textbf{Remerciements~:}~les auteurs tiennent à remercier l'équipe de la licence \textit{sciences et humanités} de l'université d'Aix-Marseille, sans qui ce travail n'aurait jamais vu le jour ; nous remercions en particulier Sara Ploquin-Donzenac pour son aide précieuse et constante. Nous tenons également à remercier  Valérie Debuiche, Christian Houzel et Sylvie Pic pour les fructueuses discussions que nous avons au sujet du \textit{Brouillon Project}. Nous remercions également les auteurs et contributeurs du logiciel Geogebra, qui nous a permis de réaliser les figures de cet article. Enfin, nous remercions plus particulièrement Philippe Abgrall pour ses nombreuses suggestions et sa relecture attentive du présent texte.

\bibliographystyle{plain}

\vs
\noindent{\bf Les auteurs~:~}\\
\noindent Marie Anglade\\
Université d'Aix-Marseille\\
CEPERC UMR CNRS 7304\\
3, place victor Hugo\\
13331 Marseille cede 3\\
France. \\
\verb=marie.anglade@univ-amu.fr=

\vs
\noindent Jean-Yves Briend\\
Université d'Aix-Marseille\\
I2M UMR CNRS 7373\\
CMI\\
39, rue Joliot-Curie\\
13453 Marseille cedex 13\\
France.\\
\verb=jean-yves.briend@univ-amu.fr=

\end{document}